\documentclass[11pt]{amsart}
\usepackage[T1]{fontenc} 
\usepackage{amsmath,amsthm,amsfonts,amssymb,bm,graphicx,xcolor,mathtools}

\usepackage
{hyperref}
\hypersetup{colorlinks=true,citecolor=blue,linkcolor=blue,urlcolor=blue,
pdfstartview=FitH }

\theoremstyle{plain}
\newtheorem{theorem}{Theorem}

\newtheorem{lemma}{Lemma}
\newtheorem{prop}{Proposition}
\mathtoolsset{showonlyrefs=true}
\theoremstyle{definition}

\renewcommand{\geq}{\geqslant}
\renewcommand{\leq}{\le}

\newcommand{\rank}{\operatorname{rank}}

\newcommand{\bdd}{\begin{center}\begin{tikzcd}}
\newcommand{\bd}{\begin{tikzcd}}
\newcommand{\edd}{\end{tikzcd}\end{center}}
\newcommand{\ed}{\end{tikzcd}}
\newcommand{\bdp}{\begin{center}\begin{tikzpicture}}
\newcommand{\edp}{\end{tikzpicture}\end{center}}
\newcommand{\bi}{\begin{itemize}}
\newcommand{\ei}{\end{itemize}}
\newcommand{\bt}{\begin{tikzpicture}}
\newcommand{\et}{\end{tikzpicture}}
\newcommand{\ba}{\[\begin{aligned}}
\newcommand{\ea}{\end{aligned}\]}
\newcommand{\bp}{\begin{pmatrix}}
\newcommand{\ep}{\end{pmatrix}}
\newcommand{\bv}{\begin{vmatrix}}
\newcommand{\ev}{\end{vmatrix}}
\newcommand{\bb}{\begin{bmatrix}}
\newcommand{\eb}{\end{bmatrix}}
\newcommand{\bB}{\begin{Bmatrix}}
\newcommand{\eB}{\end{Bmatrix}}
\newcommand{\bea}{\begin{enumerate}[label=(\alph*)]}
\newcommand{\ber}{\begin{enumerate}[label=(\roman*)]}
\newcommand{\ben}{\begin{enumerate}[label=(\arabic*)]}
\newcommand{\ee}{\end{enumerate}}

\numberwithin{equation}{section}

\makeatletter
\def\Ddots{\mathinner{\mkern1mu\raise\p@
\vbox{\kern7\p@\hbox{.}}\mkern2mu
\raise4\p@\hbox{.}\mkern2mu\raise7\p@\hbox{.}\mkern1mu}}
\makeatother
\makeatletter
\DeclareRobustCommand\widecheck[1]{{\mathpalette\@widecheck{#1}}}
\def\@widecheck#1#2{%
    \setbox\z@\hbox{\m@th$#1#2$}%
    \setbox\tw@\hbox{\m@th$#1%
       \widehat{%
          \vrule\@width\z@\@height\ht\z@
          \vrule\@height\z@\@width\wd\z@}$}%
    \dp\tw@-\ht\z@
    \@tempdima\ht\z@ \advance\@tempdima2\ht\tw@ \divide\@tempdima\thr@@
    \setbox\tw@\hbox{%
       \raise\@tempdima\hbox{\scalebox{1}[-1]{\lower\@tempdima\box
\tw@}}}%
    {\ooalign{\box\tw@ \cr \box\z@}}}
\makeatother

\begin{document}

\author{Valentin Blomer}
\address{Mathematisches Institut, Endenicher Allee 60, 53115 Bonn, Germany}
\email{blomer@math.uni-bonn.de}
 
\author{Junxian Li}
\address{Mathematisches Institut, Endenicher Allee 60, 53115 Bonn, Germany}
\email{jli135@math.uni-bonn.de}

 \title{Correlations of values of random diagonal forms}

\thanks{Both authors are supported by Germany's Excellence Strategy grant EXC-2047/1 - 390685813 as well as DFG grant BL 915/5-1}


\begin{abstract}  We study the value distribution of   diagonal forms in $k$ variables and degree $d$ with random real coefficients and positive integer variables, normalized so that the mean spacing is one. We show that the $\ell$-correlation of almost all such forms  is Poissonian when $k$ is large enough depending on $\ell$ and $d$.
\end{abstract}


\setcounter{tocdepth}{2}  \maketitle 

\maketitle

\section{Introduction}

Let  $2 \leq d \leq k$ be positive integers. For ``randomly chosen'' $\bm\alpha  \in \Bbb{R}_{>0}^k$ we consider the values at positive integer arguments of the  diagonal form 
\begin{equation}\label{form}
   q_{\bm\alpha}(\bm x) = \sum_{i=1}^k \alpha_i x_i^d
   \end{equation}
of degree $d$ in $k$ variables. If we order the non-negative values of $q_{\bm \alpha}(\bm x)^{k/d}$, $\bm x \in \Bbb{Z}_{>0}^k$, in increasing size, then the average spacing is a non-zero constant depending on $k, d, \bm \alpha$. Let $c(d, k, \bm \alpha)$ be chosen such that the mean spacing of $c(d,k, \bm \alpha)q_{\bm \alpha}(\bm x)^{k/d}$ is one. Then one would expect that ``typically'' the distribution of gaps is Poissonian, i.e.\ the probability that a randomly chosen interval $[X, X+r]$ of fixed length $r$ contains exactly $\nu$ of the numbers in the multiset $$\{\Lambda_1 \leq \Lambda_2 \leq \ldots \}  = \{c(d, k, \bm \alpha) q_{\bm \alpha}(\bm x)^{k/d} \mid \bm{x} \in \Bbb{Z}_{>0}^k\}$$ is $(\nu !)^{-1} r^{\nu} e^{-r}$. It's easy to construct counterexamples. For instance if $\bm \alpha \in \Bbb{Z}^k$ and $k > d$, then there are high multiplicities, so that the formula becomes wrong for $ r\rightarrow 0$, and if we remove the multiplicities, then the sequence is a subset of $\{n^{k/d} \mid n \in \Bbb{Z}_{>0}\}$ and again the formula becomes wrong since all the gaps are greater than one. However, for \emph{typical} ${\bm \alpha}$ (i.e.\ almost all in a Lebesgue sense, and as seen above this the best we can hope for), there is no reason to expect that this formula should fail. Unfortunately this seems to be far out of reach for any choice of parameters $2 \leq d \leq k$, but one can test it quite accurately by considering  the $\ell$-correlation function for various values of $\ell \geq 2$: for a fixed $\ell-1$-dimensional box $I = I_2 \times \ldots \times I_{\ell} \subseteq \Bbb{R}^{\ell-1}$ one would hope  to establish the asymptotic relation
\begin{equation}\label{corr}
\begin{split}
&\mathcal{T}_{\ell}(M; I; {\bm \alpha},k,d)\\
&: =  \frac{1}{M}\#\big\{(i_1, \ldots, i_{\ell}) \text{  pairwise distinct} \mid \Lambda_{i_j} - \Lambda_{i_{1}} \in I_j, \, 2 \leq j \leq \ell, \, i_j \leq M\big\} \\
&\sim \text{vol}(I) 
\end{split}
\end{equation}
as $M \rightarrow \infty$, for almost all ${\bm \alpha}$. (Note that the $\Lambda_j$ depend on $\bm \alpha$, $k$ and $d$.) 

When $d = 2$, this has an interpretation in the language of spectral geometry. In this case, the numbers $\Lambda_j$ are (rescaled) Laplace eigenvalues of a $k$-dimensional rectangular torus, and the above is precisely the celebrated Berry-Tabor conjecture \cite{BT}: for integrable systems, the local statistics of eigenvalues should come from a Poisson process, at least in generic cases. 

The first to study this type of problem was Sarnak \cite{Sa} who considered the case of pair correlation (i.e.\ $\ell = 2$) of random binary quadratic forms (i.e.\ $d=k=2$), not necessarily diagonal. A deterministic version was obtained in the breakthrough paper \cite{EMM} where pair correlation was obtained specificially for all ``diophantine'' 2-dimensional tori -- up until now this is only deterministic result of this kind.  The deterministic case of the pair correlation of the values of \emph{inhomogeneous} quadratic forms $ \sum_{i=1}^k(x_i-\alpha_i)^2$ was established by Marklof \cite{Ma1, Ma2}. The triple correlation (i.e.\ $\ell = 3$) for binary quadratic forms in a weak average sense was obtained only recently in \cite{ABR}. 

VanderKam \cite{VdK1, VdK2, VdK3} generalized Sarnak's analysis to higher correlations, higher dimensions and higher degree forms (not necessarily diagonal), but not all at the same time. The paper \cite{VdK2} deals with the pair correlation for higher degree forms, the paper \cite{VdK3} deals with higher correlations for quadratic forms. Note that M\"uller \cite{Mu} considered higher correlations of $q_{\bm \alpha}(\bm x)$ for $d=2$, which is much easier due to the different scaling unless $k=2$. (See \cite[p.238]{Mus} for a comparison.)
As an aside we remark that the investigation of small gaps in the (normalized) value distribution of (ternary) forms of arbitrary degree has been initiated by Schindler \cite{Sch}. 

\medskip

In this paper we investigate   higher correlations and restrict our attention to 
 the natural subset of \emph{diagonal}  $k$-ary forms of degree $d$ as in \eqref{form}. While superficially similar, 
 this makes the problem substantially harder for two reasons: on the one hand we have fewer coefficients to average over (namely $k$ instead of $\binom{k+d-1}{d}$),  and on a   more technical level the mixed terms can often be used to linearize the problem to some extent. Having only $d$-th powers available, the problem boils down to counting matrices with given ranks and/or subdeterminants whose entries are $d$-th powers (or differences of $d$-th powers). Unlike in previous works there is no simple transformation to understand such matrices using geometry of numbers directly. Before we talk about the methods in more detail, we state our main results.

\begin{theorem}\label{thm1} Let $d, \ell  \geq 2$ and  
\begin{equation}\label{cond}
   k > \max\Big( 2\ell ( d+ 1) - 1,B(\ell, d)\Big).
   \end{equation}
where
   \begin{align}
   	B(\ell, d)= \begin{cases} M^{ \frac{1}{4}(d + \ell+1)^2 }, & d \leq \ell - 3, \\ M^{(\ell - 1)(d+2) }, &  d > \ell - 3.\end{cases}
   \end{align}
Let $\mathcal{D} \subseteq \Bbb{R}_{>0}^k$, $I \subseteq \Bbb{R}^{\ell-1}$ be fixed boxes. 
Then there exists some $\delta>0$ such that 
\begin{equation}\label{l2}
\int_{\mathcal{D}} \big| \mathcal{T}_{\ell}(M; I; {\bm \alpha},k,d) - \text{{\rm vol}}(I) \big|^2 d\bm \alpha = O(M^{-\delta})
\end{equation}
and for almost all ${\bm \alpha}$ the $\ell$-correlation of the sequence $c(d,k, \bm \alpha)q_{\bm \alpha}(\bm x)^{k/d}$ is Poissonian. 
\end{theorem}

As one of the referees remarked, some condition on $k, \ell, d$ as in \eqref{cond} is necessary to ensure the validity of \eqref{l2}. Indeed, the $\ell$-correlation function $\mathcal{T}_{\ell}$ blows up near rational points with small denominators, so that the $L^2$ approach to almost everywhere convergence break down for $\ell$ sufficiently large (cf.\ e.g.\ \cite[Section 4]{RuS}). Quantitatively, for $\bm \alpha = (1, \ldots, 1) + O(M^{-k})$ and $k > d$ we expect $\mathcal{T}_{\ell}$ to be of size $M^{(k-d)\ell + d - 1}$, so that this set of $\bm \alpha$  can contribute as much as 
$$M^{2((k-d)\ell + d - 1) - k^2}$$
to the left hand side of \eqref{l2}. For instance, for fixed $d$, this implies a necessary condition roughly of the shape $k > 2\ell$. This is considerably weaker than \eqref{cond}, which is linear in $\ell$ while \eqref{cond} is quadratic in $\ell$. 

With more work, the number of variables in Theorem \ref{thm1} can be reduced. We demonstrate this in the case of pair correlation  where we show:
\begin{theorem}\label{thm2}  Let $\ell = 2$, $k \geq d+2$. 
Then there exists some $\delta>0$ such that \eqref{l2} holds, and for almost all ${\bm \alpha}$ the sequence $c(d, k, \bm \alpha)q_{\bm \alpha}(\bm x)^{k/d}$ has Poissonian pair correlation. 
\end{theorem}

A slightly simpler problem is to consider only a mean value rather than $L^2$-convergence. This has the benefit of requiring fewer variables. 

\begin{theorem}\label{thm3}  Let $d, \ell \geq 2$ and  
	\begin{align}\label{condweak}
	k > \max\Big( \ell ( d+ 1) - 1, B(2\ell, d)\Big).
	\end{align}
	where $B(\ell, d)$ is as in Theorem \ref{thm1}. 
Let $\mathcal{D} \subseteq \Bbb{R}_{>0}^k$, $I \subseteq \Bbb{R}^{\ell-1}$ be fixed boxes. 
Then there exists some $\delta>0$ such that 
\begin{equation}\label{weak}
\int_{\mathcal{D}} \big( \mathcal{T}_{\ell}(M; I; {\bm \alpha},k,d) - \text{{\rm vol}}(I) \big) d\bm \alpha=O(M^{-\delta}).
\end{equation}
\end{theorem}
 
Just as in Theorem \ref{thm1}, a necessary condition for $d, \ell, k$ is $ (k-d)\ell+d-1-k^2<0,$ which always holds when $k=d$. For fixed $d$ and $k>d$, this is roughly of the shape $k>\ell$, which is a weaker requirement for $k$ than that in \eqref{condweak}. But   again one can tighten the screws at several places to improve \eqref{condweak}, and we demonstrate this for the triple correlation: 
\begin{theorem}\label{thm4}  Let $\ell = 3$,   $d\geq 2,$ and $k  > \frac{4}{3}\big(d + \frac{1}{d}\big)$.  
Then there exists some $\delta>0$ such that \eqref{weak} holds.
 \end{theorem}

The relation \eqref{weak} for $\ell = 3$ together with almost everywhere pair correlation has a nice application as observed in \cite{ABR}: if a sequence with mean spacing one satisfies triple correlation in the weak sense as in Theorem \ref{thm4}, then it has infinitely many gaps of size at least $2.006$. The significance in this number lies in the fact that there exist sequences with mean spacing one and Poissonian \emph{pair} correlation with maximal gap bounded by 2. The fact that we can obtain longer gaps shows that such sequences must be genuinely closer to a Poisson process than just having Poissonian pair correlation. 

\medskip

Returning to the classical case $d=2$ in which case the considered values are rescaled eigenvalues of $k$-dimensional rectangular tori, we can further reduce the number of variables and deal with random ternary diagonal forms
$$\alpha_1 x_1^2 + \alpha_2 x_2^2 + \alpha_3 x_3^2.$$
\begin{theorem}\label{thm5}  Let $d = 2$, $k \geq 3$. Then there exists some $\delta>0$ such that \eqref{l2} holds for $\ell=2$ and \eqref{weak} holds for $\ell=3$.
 \end{theorem}
For comparison, the main result in \cite{VdK1} establishes \eqref{l2} for the forms
$$\alpha_1 x_1^2 + \ldots + \alpha_4x_4^2 + \alpha_5x_1x_2 + \ldots + \alpha_{10} x_3 x_4$$
with 10 coefficients. 

\medskip

While the investigation of forms with ``random'' real coefficients looks like a problem in classical real analysis at first sight, as soon as we consider almost-all-results and integrate over ${\bm \alpha}$ it becomes a purely arithmetic problem that needs to be tackled by tools from various branches of number theory.  Roughly speaking, we would like to replace a sum over the integer arguments $\bm{x} \in \Bbb{Z}_{>0}^k$ by an integral. This is precisely what \eqref{corr} suggests. In ``typical'' regions we can do this by the Euler McLaurin formula or Poisson summation or a similar tool. For instance, in the situation of Theorem \ref{thm3} where we have $\ell$ copies of the form \eqref{form}, this depends on the matrix
\begin{equation*}
\left(\begin{matrix} x^d_{1,1} & \cdots  &x^d_{1, k} \\ 
  \vdots & & \vdots \\ x^d_{\ell,1}  &  \cdots & x^d_{\ell, k}
\end{matrix}\right) \in \Bbb{Z}_{>0}^{\ell \times k}
\end{equation*}
or variations thereof. If this matrix has an $\ell \times \ell$ submatrix with sufficiently large determinant, we are in good shape to apply the Euler MacLaurin formula. Thus we need to show that such a matrix has ``rarely'' only small $\ell \times \ell$ subdeterminants, including the cases where it has rank strictly less than $\ell$. This is an interesting problem in its own right, but we are not aware of a systematic study of such questions. In particular, in contrast to \cite{VdK1, VdK2} such matrix counting problems do not seem amenable to a direct application of the geometry of numbers. 

 We close this introduction with the remark that \eqref{form} can also be studied with random \emph{integral} coefficients. This is a very different, but equally interesting problem, and we refer the interested reader to \cite{BD}. 
 
 \medskip
 
 \textbf{Acknowledgement:} We would like to thank the referees for a careful reading of the manuscript and useful suggestions. 

\section{The arithmetic of diagonal forms}

In this section we compile a number of results on diophantine equations and inequalities for future reference. We start with an investigation of the diophantine equation
\begin{equation}\label{eq}
a_1 x_1^d + \ldots + a_k x_k^d = 0. 
\end{equation}
with $k, d \geq 2$ and $a_1, \ldots, a_k \in \Bbb{Z} \setminus \{0\}$.

\begin{lemma}\label{lemma-a} {\rm a)} Suppose that  $k \geq 3$. Then number of integral solutions to  \eqref{eq} 
with $|x_j| \leq M$ is $O_{d, k, \varepsilon}(M^{k-2+\varepsilon})$ for any $\varepsilon>0$, uniformly in $a_1, \ldots, a_k$. 

{\rm b)} Suppose that $k=3$. Then the number   of \emph{primitive} integral solutions to \eqref{eq} with $|x_j|\leq M$ is 
 $O_{d,\varepsilon}(M^{2/d+\varepsilon})$ for any $\varepsilon > 0$, 
uniformly in the coefficients $a_1, \dots, a_3$.

{\rm c)} Suppose that $d=2$, $k=3$. Then the number   of \emph{primitive} integral solutions to \eqref{eq} with $|x_j| \leq M$ is 
$$ \ll_{\varepsilon} \Big(1 +  \frac{M\gcd(a_1a_2, a_1a_3, a_2a_3)^{1/2}}{|a_1a_2a_3|^{1/3}}\Big)|a_1a_2a_3|^{\varepsilon}$$
for any $\varepsilon > 0$. 
\end{lemma}

\textbf{Proof.} Part a) is \cite[Corollary]{BHB}, which completes the proof for all cases of $k, d$ after a long series of papers. 
 Part b) is \cite[Theorem 3]{HB}. Part c) is  \cite[Corollary 2]{BHB2}.\\

Next we study for $H, M \geq 1$, $k, d \geq 2$ and $\bm a=(a_1, \ldots, a_k) \in \Bbb{Z}^k$ the number $\mathcal{N}_{\bm{a}}(M, H)$ defined by
\begin{align}
\mathcal N_{\bm a}(M, H):=\#\Big\{(x_j)_{j=1}^k: x_j\in \mathbb{Z},\   x_j\asymp M,\  |a_1 x_1^d + \ldots + a_k x_k^d | \leq  H\Big\}
\end{align} 
with $x_j \asymp M$ by which we mean $M \leq x_j \leq 2M$.

\begin{lemma}\label{lemma-b} {\rm a)} Suppose that $(a_1, \ldots, a_k) \not= \bm 0$. Then 
$$\mathcal{N}_{\bm{a}}(M, H)\ll  M^{k-1}\Big(\frac{H}{M^{d-1}\max_i |a_i| } + 1\Big). $$

{\rm b)} If $0<|a_1|\leq |a_2|$, then
$$\mathcal N_{(a_1, -a_1, a_2, -a_2)}(M, H) \ll  M^{2+\varepsilon}  \Big(\frac{HM^{2-d}}{|a_1a_2|^{1/2}} + \frac{H^{1/2}M^{1-d/2}}{|a_1|^{1/2}} + 1\Big)$$
for all $\varepsilon> 0$. 
In the special case $d=2$, we have
$$\mathcal N_{(a_1, -a_1, a_2, -a_2)}(M,H)\ll M^{\varepsilon}\Big(\frac{M^2 H}{|a_2|}+\frac{MH}{|a_1|}+M^2\Big)$$
for all $\varepsilon>0$.

{\rm c)} If $0< |a_1|\leq |a_2|\leq |a_3|$, then 
$$\mathcal N_{(a_1,a_2,a_3)}(M, H) \ll M^{\varepsilon}\Big( \frac{M^{3-d}H}{|a_1a_3|^{1/2}}+ \frac{M^{3-d/2}H^{1/2}}{|a_3|^{1/2}}+\frac{M^{5/2-d/2}H^{1/2}}{|a_1|^{1/2}}+M^{3/2}\Big)$$
for all $\varepsilon>0$.
In the special case $d=2$, we have
\begin{align}
\mathcal N_{(a_1,a_2,a_3)}(M,H)&\ll M^\varepsilon\Big( \frac{M^2H}{|a_2|}+\frac{HM}{|a_1|}+M^2\Big)^{1/2}\Big(M+\frac{H}{|a_3|}\Big)^{1/2}
\end{align}
for all $\varepsilon>0$.
Here all the implied constants depend on $\varepsilon, k$ and $d$.
\end{lemma}

\textbf{Proof.} a)  Suppose without loss of generality that $a_1 = \max_i |a_i|$ and write $Y = |a_2 x_2^d + \ldots + a_k x_k^d|$. For fixed $x_2, \dots, x_k$ we have 
\begin{equation}\label{x1}
|x_1| = \Big(\frac{Y + O(H)}{|a_1|} \Big)^{1/d} = \Big(\frac{Y}{|a_1|}\Big)^{1/d}  + O\Bigg(\min \Big( \Big( \frac{H}{|a_1|}\Big)^{1/d}, \frac{H}{Y} \Big(\frac{Y}{|a_1|}\Big)^{1/d}\Big)+1\Bigg)
\end{equation}
where the implied constant depends only on $d$. 
If $H \gg |a_1| M^d$, then \eqref{x1} gives $$\mathcal{N}_{\textbf{a}}(M, H) \ll  M^{k-1} \Big(\frac{H}{|a_1|}\Big)^{1/d} \ll M^{k-1} \frac{H}{|a_1|M^{d-1}}.$$  If $H \ll |a_1| M^d$ with a sufficiently small implied constant, then whenever we have a solution we must have $Y  \asymp |a_1|M^d$, and \eqref{x1} yields
$$\mathcal{N}_{\textbf{a}}(M, H)\ll M^{k-1}\Big(\frac{H}{Y}\Big(\frac{Y}{|a_1|}\Big)^{1/d}+1\Big) \ll M^{k-1} \Big(\frac{H}{|a_1| M^{d-1}} + 1\Big),$$
which completes the proof.

b) For  $h \in \Bbb{Z}$ define $n_{a,b}(M,h)$ by
\begin{align}
n_{a,b}(M,h):=\#\Big\{ (x_i)_{i=1}^4: x_j\in \mathbb{Z},\ x_i\asymp M, \ a (x_1^d - x_2^d) + b (x_3^d - x_4^d) = h\Big \}.
\end{align}
Using the F\'ejer kernel  and the Cauchy-Schwarz inequality we have 
\begin{displaymath}
	\begin{split}
		&\mathcal{N}_{(a_1, -a_1, a_2, -a_2)}(M, H)  \leq 2\sum_{|h| \leq 2H}  \frac{2H - |h|}{2H}  n_{a_1, a_2}(M, h)  \\
		& = 2\sum_{|h| \leq 2H}  \frac{2H - |h|}{2H}  \sum_{x, y, z, w \asymp M}  \int_{0}^1e((a_1 (x^d - y^d) + a_2 (z^d - w^d) - h)\beta) d\beta\\
		& = 2 \int_{0}^1  \frac{1}{2H} \Big( \frac{\sin(\pi 2 H \beta)}{\sin (\pi \beta)}\Big)^2 \Big| \sum_{x \asymp M}  e(a_1x^d \beta) \Big|^2 \Big| \sum_{x \asymp M}  e(a_2x^d \beta) \Big|^2  d\beta\\
		& \leq \Big(2 \sum_{|h| \leq 2H/|a_1|} \frac{2H - |a_1h|}{2H}  n_{1, 1}(M,h) \Big)^{1/2} \Big(2 \sum_{|h| \leq 2H/|a_2|} \frac{2H - |a_2h|}{2H}  n_{1, 1}(M,h) \Big)^{1/2} \\
		& \leq \Big(2 \mathcal{N}_{(1, -1, 1,-1)}(M, 2H/|a_1|)\Big)^{1/2} \Big(2 \mathcal{N}_{(1, -1, 1,-1)}(M, 2H/|a_2|)\Big)^{1/2}. 
	\end{split}
	\end{displaymath}
		The lemma follows from \cite[Theorem 2]{RoS}, which implies 
		\begin{align}\label{RSthm}
		\mathcal N_{(1,-1, 1-1)}(M, \delta M^d)\ll_\varepsilon M^{2+\varepsilon}+ \delta M^{4+\varepsilon}
		\end{align} for any $\delta>0$ and any $\varepsilon>0$.\\
	
		When $d=2$, the result follows from a divisor argument. More precisely, we can fix $x_1, x_2$ with $x_1\not=x_2$ to obtain $O(M^\varepsilon H/|a_2|)$ choices for $x_3, x_4$ with $x_3\not=x_4$. If $x_1\not=x_2$ and $x_3=x_4$, then we have $O(M^\varepsilon H/|a_1|)$ choices for $x_1,x_2$. Finally we have $O(M^2)$ choices with $x_1=x_2$ and $x_3=x_4$. \\
	
c) By the same argument as in part b)   with H\"older's inequality, we have 
\begin{displaymath}
\begin{split}
\mathcal N_{(a_1,a_2,a_3)}(M, H)&\ll \mathcal N_{(a_1,-a_1, a_1, -a_1)}(M,2H)^{1/4} \mathcal N_{(a_2,-a_2, a_2, -a_2)}(M,2H)^{1/4}\mathcal N_{(a_3, -a_3)}(M, 2H)^{1/2}\\
& \ll \mathcal N_{(1, -1, 1, -1)}(M,2H/|a_1|)^{1/2}  \mathcal N_{(1, -1)}(M,2H/|a_3|)^{1/2}. 
\end{split}
\end{displaymath}
Clearly 
\begin{align}\label{N11}
\mathcal{N}_{(1, -1)}(M, H) \ll \#\{x, y \asymp M :  |x-y| \ll M^{1-d}H \} \ll M +  HM^{2-d},
\end{align} and the first claim then follows again from \eqref{RSthm}.  

If $d=2$, we can use the bound in part b) after an application of the Cauchy-Schwarz inequality. More precisely, we have
\begin{align}
\mathcal N_{(a_1,a_2,a_3)}(M, H)&\ll \mathcal N_{(a_1,-a_1, a_2, -a_2)}(M,2H)^{1/2}\mathcal N_{(a_3, -a_3)}(M,2H)^{1/2}, \\
&\ll \mathcal N_{(a_1,-a_1, a_2, -a_2)}(M,2H)^{1/2}\mathcal N_{(1, -1)}(M,2H/|a_3|)^{1/2}, 
\end{align}
and we can apply the case $d=2$ in part b) and \eqref{N11} to obtain the second claim.

\section{Reduction of the problem}\label{sec3}

As mentioned in the introduction, at the heart of the problem is a number theoretic question on the behaviour of integral matrices with $d$-th power entries. In this section we reduce \eqref{l2} and \eqref{weak} to such arithmetic problems. This reduction is fairly standard (see e.g.\ \cite{Sa, VdK1, VdK2, VdK3}), and we can be brief. From now on all implied constants may depend on $d, \ell, k, \varepsilon$ where applicable. 

For the Fourier and diophantine analysis to come we need to smooth out and slightly simplify the expression $\mathcal{T}_{\ell}(M; I; {\bm \alpha},k,d)$. We replace the conditions  
$$\Lambda_{i_j} - \Lambda_{i_1} = c(d, k, {\bm \alpha}) \big(q_{\bm \alpha}(\bm{x}_j)^{k/d} - q_{\bm \alpha}(\bm{x}_1)^{k/d}\big) \in I_j, \quad 2 \leq j \leq \ell.$$ 
by smooth and slightly rescaled versions
$$ W_j\big((q _{\bm\alpha}(\bm x_j) - q _{\bm\alpha}(\bm{x}_1))M^{k-d}\big) $$
where the $W_j$'s are fixed smooth, compactly supported weight functions. We also smooth out the summation over the variables $\bm{x}_1, \ldots, \bm{x}_{\ell} \in \Bbb{Z}_{>0}^k$ using  smooth weight functions $\Psi_j$, $1 \leq j \leq \ell$, that in addition have compact support in $[1, 2]$, so that all the variables are localized. With this in mind we define the quantity
$$\mathcal{T}^{\ast}_{\ell}(M; {\bm \alpha},k,d):= \underset{\bm x_1, \ldots, \bm x_{\ell} \in \Bbb{Z}_{>0}^k}{\left.\sum\right.^{\ast}}    \prod_{j=2}^{\ell} W_j\big((q _{\bm\alpha}(\bm x_j) - q _{\bm\alpha}(\bm{x}_1))M^{k-d}\big) \prod_{j=1}^{\ell} \Psi_j\Big(\frac{\bm x_j}{M}\Big) $$
where the asterisk indicates that the sum is over pairwise distinct vectors $\bm x_j = (x_{j, i}) \in \Bbb{Z}_{>0}^k$. This should then be compared with the (``Hardy-Littlewood'') expectation
$$M^{k} \Big(\prod_{j=2}^{\ell} \widehat{W}_j(0)\Big) c(\bm \alpha)$$
where $\widehat{W}_j$ denotes the Fourier transform and 
 $$c(\bm \alpha) = \int_{\Bbb{R}^{\ell-1}} \int_{\Bbb{R}^{k\ell}}e\Big(\sum_{j=2}^{\ell} \xi_j\big(q_{\bm \alpha}(\bm{x}_j) - q_{\bm \alpha}(\bm{x}_1)\big)\Big) \prod_{j=1}^{\ell} \Psi_j(\bm{x}_j) d\bm{x}\, d \bm \xi.$$
This constant can also be interpreted as the surface integral
$$\int_{\mathcal{Q}(\bm \alpha)} \prod_{j=1}^{\ell} \Psi_j(\bm{x}_j) d\bm{x}$$
 where $\mathcal{Q}(\bm \alpha)$ is the surface given by $q_{\bm \alpha}(\bm{x}_j) - q_{\bm \alpha}(\bm{x}_1) = 0$, $2 \leq j \leq \ell$. Finally we also smooth out the region $\mathcal{D}$ over which we integrate $\bm \alpha$ using a smooth test function $F$ with compact support\footnote{We think of $\bm \alpha$ as a column  vector, so $F$ takes column vectors as arguments. Sometimes we write $\bm \alpha = (\begin{smallmatrix} \bm \alpha_1 \\ \bm \alpha_2 \end{smallmatrix})$ with two column vectors of lengths $\ell$ and $k-\ell$ respectively, in which case we write $F(\bm \alpha_1, \bm \alpha_2)$ with the obvious meaning.}
  on $\Bbb{R}_{>0}^k$. 
 
 By standard analytic methods (see \cite[Section 3]{Sa} or \cite[Section 2.3]{VdK2} and also \cite{VdK1, VdK3})  one shows the following.
 \begin{prop}\label{prop1} Let $\Psi_j$, $W_j$, $F$ as described above, $M \geq 1$, and fix a triple $(d, \ell, k)$ with $2 \leq d, \ell \leq k$.  
 
 {\rm a)} Suppose there exists some $\delta > 0$ such that 
 \begin{equation}\label{l2simple}
 \int_{\Bbb{R}_{>0}^k} F(\bm \alpha)\Big[\mathcal{T}^{\ast}_{\ell}(M; {\bm \alpha},k,d) - M^k\Big(\prod_{j=2}^{\ell} \widehat{W}_j(0) \Big)  c(\bm \alpha)\Big]^2 d\bm \alpha \ll M^{2k - \delta}.
 \end{equation}
 Then \eqref{l2} holds   and for almost all ${\bm \alpha}$ the $\ell$-correlation of the sequence $c(d,k, \bm \alpha)q_{\bm \alpha}(\bm x)^{k/d}$ is Poissonian. 
 
{\rm  b)} Suppose there exists some $\delta > 0$ such that 
 \begin{equation}\label{weaksimple}
 \int_{\Bbb{R}_{>0}^k} F(\bm \alpha)\Big[\mathcal{T}^{\ast}_{\ell}(M; {\bm \alpha},k,d) - M^k\Big(\prod_{j=2}^{\ell} \widehat{W}_j(0) \Big)  c(\bm \alpha)\Big] d\bm \alpha \ll M^{k - \delta}.
 \end{equation}
 Then \eqref{weak} holds.
 \end{prop}

To analyze \eqref{weaksimple} we write
\begin{equation}\label{quantity}
\begin{split}
 \mathcal{C}(M) &:= \mathcal{C}(M; k, \ell, d)  :=   \int_{\Bbb{R}_{>0}^k} F(\bm \alpha)   \mathcal{T}^{\ast}_{\ell}(M; \bm \alpha, k, d) d\bm \alpha \\
 &=  \underset{\bm x_1, \ldots, \bm x_{\ell} \in \Bbb{Z}_{>0}^k}{\left.\sum\right.^{\ast}}   \int_{\Bbb{R}_{>0}^k} F(\bm \alpha)  \prod_{j=2}^{\ell} W_j\big((q _{\bm\alpha}(\bm x_j) - q _{\bm\alpha}(\bm{x}_1))M^{k-d}\big) \prod_{j=1}^{\ell} \Psi_j\Big(\frac{\bm x_j}{M}\Big) d\bm \alpha\\
& =: \underset{\bm x_1, \ldots, \bm x_{\ell} \in \Bbb{Z}_{>0}^k}{\left.\sum\right.^{\ast}}  
 \mathcal I(\bm x, M),
\end{split}
\end{equation}
say (where of course also $\mathcal{I}$ depends on $k, \ell, d$). We will usually drop the dependence on $k, \ell, d$ from the notation. We  define for $\bm x=(\bm x_j)\in \Bbb{Z}_{>0}^{\ell\times k}$ the matrix
\begin{equation}\label{TT'}
T = T(\bm x)= \left(\begin{matrix} x^d_{1,1} & \cdots  &x^d_{1 ,k} \\ x^d_{2,1} - x_{1,1}^d  &  \cdots & x^d_{2,k} - x_{1,k}^d
\\  \vdots & & \vdots \\ x^d_{\ell,1} - x_{1,1}^d  &  \cdots & x^d_{\ell, k} - x_{1,k}^d
 \end{matrix}\right) \in \Bbb{Z}^{\ell \times k}
 \end{equation}
 and write $T = (T_1\, T_2) \in \Bbb{Z}^{\ell \times (\ell + (k-\ell))}$ so that $$T_1 = T_1(\bm{x})$$ is a square matrix. We split the sum over the $\bm x_j$   into pieces according to $r = \rank(T)$ and call the corresponding piece 
 \begin{align}\label{crdef}
 \mathcal{C}_{r}(M) = \mathcal{C}_r(M; k, \ell, d):=\underset{\substack{\bm x_1, \ldots, \bm x_{\ell} \in \Bbb{Z}_{>0}^k\\ \rank T=r}}{\left.\sum\right.^{\ast}}  \mathcal{I}(\bm{x}, M).
 \end{align}
 When $T$ is of full rank, $\mathcal{C}_{ \ell}(M)$ is the contribution of those $\bm{x}_1, \ldots, \bm{x}_\ell$ where $T_1$ is invertible (over $\Bbb{Q}$).   In this case we write $\bm\alpha = (\begin{smallmatrix} \bm \alpha_1\\ \bm \alpha_2 \end{smallmatrix}) \in \Bbb{R}^{\ell + (k-\ell)}$, and we introduce new variables $\bm a = T\bm \alpha =   T_1\bm \alpha_1    + T_2 \bm\alpha_2     = (a_1, \ldots, a_{\ell})^{\top} \in \Bbb{R}^{\ell}$ and obtain
  \begin{equation*}
  \begin{split}
  \mathcal{I}(\bm{x}, M) =  \frac{1}{|\det T_1|} \prod_{j=1}^{\ell} \Psi_j\Big(\frac{\bm x_j}{M}\Big)  \int_{\Bbb{R}^{\ell}} \int_{\Bbb{R}^{k-\ell}} F(T_1^{-1}({\bm a} - T_2 \bm \alpha_2), {\bm \alpha_2}) \prod_{j=2}^{\ell} W_j\big(a_jM^{k-d}\big)d\bm\alpha_2\, d{\bm a}.
\end{split} 
\end{equation*}
 Note that we may have some choice to decompose the matrix $T$ and extract an  $\ell \times \ell$ invertible submatrix $T_1$.  In practice  it will often be convenient to relabel the columns such that $|\det T_1|$ is maximal among all $\ell \times \ell$ minors. 
Let $D_0 = M^{d\ell - \eta}$ for some $\eta > 0$ and cut off smoothly the portion $|\det T_1 |\leq D_0$ from the sum, i.e.\ we insert a weight function $\phi(\det T_1/D_0)$  where $\phi$ is a smooth weight function with support on $(-\infty, -1]\cup[1, \infty)$ 
and define
\begin{equation}\begin{split}\label{CD+}
 \mathcal{C}_{\ell,+}(M, D_0)&=\underset{\substack{\bm x_1, \ldots, \bm x_{\ell} \in \Bbb{Z}_{>0}^k\\ \rank T=\ell}}{\left.\sum\right.^{\ast}}  \mathcal{I}(\bm{x}, M)\phi\Big(\frac{\det T_1}{D_0}\Big),\\
 \mathcal{C}_{\ell, -}(M, D_0)&= \mathcal{C}_{\ell}(M) - \mathcal{C}_{\ell, +}(M, D_0).
\end{split}  
\end{equation}

  The key observation is that  for $|\det T_1| \geq D_0$ the integral $\mathcal{I}(\bm{x}, M)$    is a ``flat'' function in each variable $x_{j,i}$, and hence by the Euler-MacLaurin summation formula (or a similar device) we can replace the sum by an integral, up to a small error term. Note that for $\rank(T) = \ell$ the extra summation conditions imposed by the asterisk are void. More precisely,  as the integral over $\bm \alpha_2$ is $O(1)$, the integral over $ a_2, \ldots, a_{\ell}$ is $O(M^{(d-k)(\ell-1)})$ and the integral over $a_1$ is $O(M^d)$, we have 
$$\| \nabla\mathcal{I}(\bm{x}, M)\| \ll \frac{1}{|\det T_1|}  M^{(\ell - 1)(d-k) + d}  \cdot \frac{M^{d\ell - 1}}{|\det T_1|} \leq \frac{M^{k-1 + \ell(d-k) + \eta}}{|\det T_1|}.$$
The $L^1$-norm of this expression over the region $x_{j,i} \asymp M$ is $M^{k-1+\eta + \varepsilon}$, so that
$$\mathcal{C}_{\ell, +}(M, D_0) = \mathcal{G}_+(M, D_0) + O(M^{k - 1 + \eta + \varepsilon})$$ 
where  $\mathcal{G}_+(M, D_0)$ is given by replacing the sum by an integral: 
$$ \int_{\Bbb{R}^{k\ell} } \frac{\phi(\det T_1/D_0) }{|\det T_1|} \prod_{j=1}^{\ell} \Psi_j\Big(\frac{\bm x_j}{M}\Big)  \int_{\Bbb{R}^{\ell}} \int_{\Bbb{R}^{k-\ell}} F(T_1^{-1}({\bm a} - T_2 \bm \alpha_2), {\bm \alpha_2}) \prod_{j=2}^{\ell} W_j\big(a_jM^{k-d}\big)d\bm\alpha_2\, d{\bm a} \, d\bm{x}.$$

 The manipulation of $\mathcal{G}_+(M, D_0)$  is now standard. By Cramer's rule we have
$$F(T_1^{-1}({\bm a} - T_2 \bm \alpha_2), {\bm \alpha_2})  - F(T_1^{-1}((a_1, 0, \ldots, 0)^{\top} - T_2 \bm \alpha_2), {\bm \alpha_2)} \ll \frac{M^{d(\ell - 1)}}{|\det T_1|} \cdot \frac{1}{M^{k-d}} \ll M^{\eta - k}$$
on the support of $W_2 \cdots W_{\ell}$, so that up to a total error of $O(M^{2\eta})$ after integration over $\bm \alpha_2$ and $ \bm a$ we may replace the first term on the left hand side of the previous display by the second in the definition of  $\mathcal{G}_+(M, D_0)$. 
We can now integrate over $a_2, \ldots, a_\ell$ and rescale the $x_{j,i}$-variables and the $a_1$-variable to obtain that  $$ \mathcal{G}_+(M, D_0) = M^k \Big(\prod_{j=2}^{\ell} \widehat{W}_j(0) \Big)\mathcal{G}_+(\eta)+ O(M^{2\eta})$$ where
$$\mathcal{G}_+ (\eta)  =
 \int_{\Bbb{R}^{k\ell} } \frac{\phi(M^{\eta} \det T_1)}{|\det T_1|} \prod_{j=1}^{\ell} \Psi_j ( \bm x_j )  \int_{\Bbb{R}} \int_{\Bbb{R}^{k-\ell}} F(T_1^{-1}((a_1, 0, \ldots, 0)^{\top} - T_2 \bm \alpha_2), {\bm \alpha_2})  d\bm\alpha_2\, da_1 \, d\bm{x}.$$
 Correspondingly we also define $\mathcal{G}_-(\eta)$ where the factor $\phi(M^{\eta} \det T_1)$ is replaced with $1-\phi(M^{\eta} \det T_1)$. 
 By  Fourier inversion (see \cite[p.\ 191]{Sa} or \cite[p.\ 389]{VdK2}) one shows that
 $$\int_{\Bbb{R}^k} F(\bm \alpha) c(\bm \alpha) d\bm \alpha = \mathcal{G}_{+}(\eta) + \mathcal{G}_-(\eta),$$
and we conclude the following.
\begin{prop}\label{prop2}
Let $k\geq \ell$. Suppose that there exist $0< \delta', \eta < 1$ and some choice of $T_1$ such that for $D_0 = M^{d\ell - \eta}$ we have with the above notation
\begin{equation}\label{prop2-bound}
\sum_{r \leq \ell - 1} \mathcal{C}_r(M) +  \mathcal{C}_{\ell, -}(M, D_0) + M^k \mathcal{G}_-(\eta) \ll M^{k - \delta'},
\end{equation}
then \eqref{weaksimple} holds for some $\delta > 0$. 
\end{prop}

The treatment of \eqref{l2simple} is very similar. After expanding the square, we will need for the mixed terms an asymptotic formula for $\mathcal{C}(M; k ,\ell, d)$ with $F(\bm \alpha)$ replaced by $F(\bm \alpha) c(\bm \alpha)$, which can be studied with Proposition \ref{prop2}. It remains to consider the critical quantity 
\begin{equation}\label{quantity2}
\begin{split}
\mathcal{C}^{(2)}(M) := & \mathcal{C}^{(2)}(M; k, \ell, d) :=  \underset{\substack{\bm x_1, \ldots, \bm x_{\ell} \in \Bbb{Z}_{>0}^k\\ \bm y_1, \ldots, \bm y_{\ell} \in \Bbb{Z}_{>0}^k}}{\left.\sum\right.^{\ast}} \int_{\Bbb{R}_{>0}^k} F(\bm \alpha) \prod_{j=1}^{\ell} \Psi_j\Big(\frac{\bm x_j}{M}\Big) \Psi_j\Big(\frac{\bm y_j}{M}\Big)  \\
& \quad\quad \times \prod_{j=2}^{\ell} W_j\big((q _{\bm\alpha}(\bm x_j) - q _{\alpha}(\bm x_1))M^{k-d}\big) W_j\big((q _{\bm\alpha}(\bm y_j) - q _{\bm\alpha}(\bm y_1))M^{k-d}\big) d\bm\alpha\\
=:&  \underset{\substack{\bm x_1, \ldots, \bm x_{\ell} \in \Bbb{Z}_{>0}^k\\ \bm y_1, \ldots, \bm y_{\ell} \in \Bbb{Z}_{>0}^k}}{\left.\sum\right.^{\ast}}\mathcal{I}^{(2)}(\bm{x}, \bm{y}, M), 
\end{split}
\end{equation}
say, where the asterisk indicates that  vectors $\bm x_j = (x_{j, i}) \in \Bbb{Z}_{>0}^k$ are pairwise distinct and the vectors $\bm y_j = (y_{j, i}) \in \Bbb{Z}_{>0}^k$ are pairwise distinct. We treat this quantity analogously. We define the matrix
\begin{equation}\label{Txy}
{\tt T} = \left(\begin{matrix} T(\bm x)\\ T(\bm y)\end{matrix}\right)  \in \Bbb{Z}^{2\ell \times k}
 \end{equation}
with $T$ as in \eqref{TT'}, and decompose $\tt T = (\tt T_1\, \tt T_2)$ with $\tt T_1 \in \Bbb{Z}^{2\ell \times 2\ell}$ a square matrix (again we have the option to  relabel the rows and columns). For $r \leq 2\ell$  we define as above $\mathcal{C}_r^{(2)}(M)$ to be the contribution of those values of $\bm{x}_1, \ldots \bm{x}_\ell, \bm{y}_1, \ldots, \bm{y}_\ell $ where $\rank({\tt T}) = r$ for $1\leq r\leq 2\ell$, i.e. 
\begin{align}\label{cr2def}
\mathcal{C}_{r}^{(2)}(M) = \mathcal{C}_r^{(2)}(M; k, \ell, d):=\underset{\substack{\bm x_1, \ldots, \bm x_{\ell} \in \Bbb{Z}_{>0}^k\\ \bm y_1, \ldots, \bm y_{\ell} \in \Bbb{Z}_{>0}^k\\\rank \tt{T}=r}}{\left.\sum\right.^{\ast}}  \mathcal{I}^{(2)}(\bm{x},\bm y, M)
\end{align}
where the asterisk has the same meaning as in \eqref{quantity2}.
For $D_0 = M^{2d\ell - \eta}$ and $\phi$ a smooth weight function with support on $(-\infty, -1]\cup[1, \infty)$, we define 
\begin{equation}\begin{split}\label{CD2+}
\mathcal{C}_{2\ell,+}^{(2)}(M, D_0)&=\underset{\substack{\bm x_1, \ldots, \bm x_{\ell} \in \Bbb{Z}_{>0}^k\\\bm y_1, \ldots, \bm y_{\ell} \in \Bbb{Z}_{>0}^k\\ \rank \tt{T}=2\ell}}{\left.\sum\right.^{\ast}}  \mathcal{I}^{(2)}(\bm{x},\bm y, M)\phi\Big(\frac{\det \tt{T}_1}{D_0}\Big),\\
\mathcal{C}_{2\ell, -}^{(2)}(M, D_0)&= \mathcal{C}_{2\ell}^{(2)}(M) - \mathcal{C}_{2\ell, +}^{(2)}(M, D_0).
\end{split}  
\end{equation}
  We can extract the expected main term by replacing  the sum over $\bm{x}_1, \ldots \bm{x}_k, \bm{y}_1, \ldots, \bm{y}_k$ by an integral featuring the constant
\begin{displaymath}
\begin{split}
\mathcal{G}^{(2)}   = &
 \int_{\Bbb{R}^{k\ell} } \int_{\Bbb{R}^{k\ell} }  \frac{1}{|\det \tt T_1|} \prod_{j=1}^{\ell} \Psi_j ( \bm x_j ) \Psi_j ( \bm y_j )\\
 &  \int_{\Bbb{R}^2} \int_{\Bbb{R}^{k-2\ell}} F( {\tt T}_1^{-1}((a_1, 0, \ldots, 0, a_{\ell+1}, 0, \ldots, 0)^\top - {\tt T}_2 \bm \alpha_2), {\bm \alpha_2}) )d\bm\alpha_2\, d(a_1, a_{\ell}) \, d\bm{x}\, d\bm{y}.
 \end{split}
 \end{displaymath}
As before we define $\mathcal{G}^{(2)}_+(\eta)$ by inserting a smooth weight $\phi(M^{\eta}\det \tt T_1)$ and $\mathcal{G}^{(2)}_-(\eta)$ by inserting a smooth weight $1-\phi(M^{\eta}\det \tt T_1)$. We also have 
\begin{align}
\int_{\mathbb{R}^k}F(\bm \alpha)c(\bm \alpha)^2d\bm \alpha=\mathcal G^{(2)}(\eta)=\mathcal G_+^{(2)}(\eta)+\mathcal G_{-}^{(2)}(\eta).
\end{align} 
Thus we obtain:
\begin{prop}\label{prop3} Let $k\geq 2\ell$. 
Suppose that there exist $0< \delta', \eta < 1$ and some choice of ${\tt T}_1$ such that for $D_0 = M^{2d\ell - \eta}$ we have 
\begin{equation}\label{prop3-bound}
\sum_{r \leq 2\ell - 1} \mathcal{C}^{(2)}_r(M) +  \mathcal{C}^{(2)}_{2\ell, -}(M, D_0) + M^{2k} \mathcal{G}^{(2)}_-(\eta) \ll M^{2k - \delta'},
\end{equation} 
and in addition \eqref{prop2-bound} holds with $D_0 = M^{d\ell - \eta}$, 
then \eqref{l2simple} holds for some $\delta > 0$. 
\end{prop}

In order to prove Theorems \ref{thm1} -- \ref{thm5}, we will verify the relevant bounds in Propositions \ref{prop2} and \ref{prop3}. To this end we need to bound the integrals $\mathcal{I}(\bm{x}, M)$ resp.\ $\mathcal{I}^{(2)}(\bm{x}, \bm{y}, M)$ (this is fairly simple) and then count the number of matrices $T$ resp.\ ${\tt T}$ satisfying certain determinant and rank conditions. This is the main arithmetic work and essentially occupies the rest of the paper.

\section{Proof of Theorem \ref{thm3}}

Before we start with the proof, we define the auxiliary matrix 
\begin{equation}\label{T'}
 T' = T'(\bm x) = \left(\begin{matrix} x^d_{1,1} & \cdots  &x^d_{1, k} \\ x^d_{2,1}    &  \cdots & x^d_{2,k} 
\\  \vdots & & \vdots \\ x^d_{\ell,1}    &  \cdots & x^d_{\ell, k} 
\end{matrix}\right) \in \Bbb{Z}_{>0}^{\ell \times k},
\end{equation}
which clearly satisfies $\rank T'=\rank T = r$ with $T$ as in \eqref{TT'}. 
 Note that we must have $r \geq 2$. In fact, more strongly we can state that any $\bm x_j$, $\bm x_{j'}$ with $j \not= j'$ cannot be linearly dependent:  they cannot be identical, and if they are different, but linearly dependent, then $q_{\bm \alpha}(\bm x_j) - q_{\bm\alpha}(\bm x_{j'}) \gg M^{d-1}\geq M$ for $\bm \alpha \gg 1$ and so they are not in the support of $W_j$ when $k\geq 2$. 

\subsection{Small ranks} We begin by studying the case when $r = \text{rank}(T) < \ell$  and give an upper bound for  $\mathcal{C}_{r}(M) = \mathcal{C}_{r}(M; k, \ell, d)$ which we recall is defined in \eqref{crdef}.

\begin{lemma}\label{lem2} Let $k\geq d, \ell \geq 2$. If  $2 \leq r < \ell$, then
$$\mathcal{C}_{ r}(M; k, \ell, d) 
\ll \begin{cases} M^{ \frac{1}{4}(d + \ell+1)^2 + \varepsilon}, & d \leq \ell - 3, \\ M^{(\ell - 1)(d+2) + \varepsilon}, &  d > \ell - 3.\end{cases}$$
for any $\varepsilon > 0$. 
\end{lemma}

\textbf{Proof.}  
Exchanging rows of $T'$ if necessary,  
we can assume that  each of the $\ell - r$ bottom rows of $T'$ is a linear combination of the first $r$ rows, i.e.
\begin{equation}\label{linear}
x^d_{j,i} = \sum_{\nu = 1}^r \rho_{j, \nu}  x_{\nu, i}^d, \quad r+1 \leq j \leq \ell, \quad 1 \leq i \leq k
\end{equation}
for suitable $\rho_{j, \nu} \in \Bbb{Q}$ (which are independent of $i$), and the first $r$ rows of $T'$ (and hence of $T$) have rank $r$. 
We note already at this point that 
\begin{equation}\label{two}
\text{at most two of the $\rho_{j, \nu}$, $1 \leq \nu \leq r$, are nonzero},
\end{equation}
for otherwise $\bm x_j$ would be a multiple of some $\bm x_{\nu}$, $1 \leq \nu \leq r$, which is impossible.

We write
 \begin{equation}\label{decomp} 
  T = \left( \begin{matrix} T_1 & T_2\\ T_3& T_4 \end{matrix}\right) \in \Bbb{Z}^{(r + (\ell-r)) \times (r + (k-r))}, 
  \end{equation}
where  $T_1 \in \Bbb{R}^{r \times r}$ is invertible (since the upper left $r\times r$ block of $T'$ is invertible).

We write $\bm \alpha =( \begin{smallmatrix} \bm\alpha_1 \\ \bm\alpha_2\end{smallmatrix}) \in \Bbb{R}^{r + (k-r)}$ in  \eqref{quantity}, and we introduce new variables by writing $T_1 \bm \alpha_1  + T_2\bm\alpha_2 =\bm a  \in \Bbb{R}^r$, so that $\bm \alpha_1 =  T_1^{-1}(\bm a  - T_2\bm \alpha_2)$.   
We re-write $\mathcal{I}(\bm{x}, M)$  as an integral over $(\begin{smallmatrix}\bm a \\ \bm \alpha_2\end{smallmatrix}) \in \Bbb{R}^{r + (k-r)}$. The integral over $\bm \alpha_2$ is $O(1)$, the integral over $a_2, \ldots, a_r$ is $O(M^{(r-1)(d-k)})$ and the remaining integral over $a_1 = q_{\bm\alpha}(\bm x_1)$ is  $O(M^d)$. 
We conclude that 
\begin{equation}\label{inte} \mathcal{I}(\bm{x}, M) \ll  M^{(r-1)(d-k) + d} |\det T_1 |^{-1}
\end{equation}
if $\text{rank}(T) = r$. 
 At this point we use the trivial bound $|\det T_1|^{-1}\ll 1$ and summarize our previous discussion by stating that 
  \begin{equation}\label{clr}
  \mathcal{C}_{r}(M ) \ll   M^{(r-1)(d-k)+ d} \mathcal{C}^{\ast}_{r}(M)
  \end{equation}
   where 
   \begin{align}
   \mathcal C^{\ast}_r(M)=\#\Big\{ (x_{j,i})_{\ell \times k}:x_{j,i}\in \mathbb{Z}, \ x_{j,i}\ll M  \text{ and there exist } \rho_{j, \nu} \text{ such that   \eqref{linear}  holds } \Big\}.
   \end{align}
    We proceed to estimate $\mathcal{C}^{\ast}_{r}(M)$. 
Suppose that the first $r$ columns of $T'$ have already been fixed. We use this in combination with \eqref{linear} to determine  the numbers $\rho_{j, \nu}$. Indeed, for each $r+1 \leq j \leq \ell$ we must have
\begin{equation*}
(T_1')^{\top} \left(\begin{matrix} \rho_{j,1}\\ \vdots \\ \rho_{j,r}\end{matrix}\right) =  \left(\begin{matrix} x_{j,1}^d  \\ \vdots \\ x_{j, r}^d  \end{matrix}\right). 
\end{equation*}
 Since $T_1$ in invertible, this determines  each $\rho_{j, \nu} \in \Bbb{Q}$.   Now we return to 
 \eqref{linear} and recall \eqref{two}.  Thus we see that  \eqref{linear} is a non-degenerate form of degree $d$ in at least 3 variables. By 
Lemma \ref{lemma-a}a, 
it follows that the number of choices for the $(r+1)$-tuple $x_{1, i}, \ldots, x_{r, i}, x_{j, i}$  is at most 
\begin{equation}\label{r-1}
   O(M^{r-1+\varepsilon}),
 \end{equation}   uniformly in the coefficients $ \rho_{j, \nu}$. Once these are chosen for some fixed $j$, e.g.\ $j = r+1$, then by \eqref{linear}  the remaining $x_{j,i}$'s are determined. We conclude that
 \begin{equation*}
\mathcal{C}^{\ast}_{ r}(M) \ll  M^{\ell r + (r-1) (k - r) + \varepsilon},
\end{equation*}
so that
$$\mathcal{C}_{r}(M) \ll M^{(r-1)(d-k) +d+ \ell r + (r-1) (k - r ) + \varepsilon} = M^{\ell r + d r + r - r^2  +\varepsilon}.$$
The worst case is $r = \min(  \frac{1}{2}(\ell + d+1), \ell - 1)$, and we obtain the lemma.

\subsection{Full rank}\label{sec-regular}

We proceed to estimate $\mathcal{C}_{\ell,-}(M, D_0)$ (which also depends on $k$ and $d$) as defined in \eqref{CD+}.

We start by arranging the columns of $T'$ defined in \eqref{T'} as follows. For any two subsets $\mathfrak{c} \subseteq \{1, \ldots, k\}$ and $\mathfrak{r} \subseteq \{1, \ldots, \ell\}$ of   cardinality $n$ we denote by $\mathcal{T}^{\mathfrak{c}}_{\mathfrak{r}}$ the determinant of the submatrix of $T'$ consisting of the columns in $\mathfrak{c}$ and rows in $\mathfrak{r}$. We may then assume that $T'$ is ``optimally arranged'' (this phrase was coined by VanderKam \cite{VdK3})  in the following sense: for each $1 \leq n \leq \ell$ we have 
\begin{equation}\label{oa}
|\mathcal{T}^{\{1, \ldots, n\}}_{\{1, \ldots, n\}} | \geq  |\mathcal{T}^{\{1, \ldots, n-1, m\}}_{\{1, \ldots, n\}} |
\end{equation}
for all $n < m \leq k$. (Then the left hand side is automatically non-zero since $\rank(T') = \ell$.) Indeed, this can be arranged by permuting columns as follows: first we move the largest entry of the first row into the upper left corner. Then we permute the columns $2, 3, \ldots, k$ so that \eqref{oa} holds for $n = 2, 3$ etc.

We return to the estimation of \eqref{quantity}. As in \eqref{decomp}, we write  $T = (T_1 \  T_2)\in \mathbb{R}^{\ell \times (\ell+(k-\ell))}$ with $\rank (T_1)=\ell$.
Again we write $\bm\alpha = (\begin{smallmatrix} \bm\alpha_1 \\  \bm\alpha_2 \end{smallmatrix}) \in \Bbb{R}^{\ell + (k-\ell)}$, and we introduce new variables $\bm a   = T_1\bm \alpha_1  + T_2\bm\alpha_2   \in \Bbb{R}^{\ell}$, so that \begin{equation}\label{quantity-new}
\begin{split}
   \mathcal{I}(\bm{x}, M) 
   &\ll M^{(\ell - 1)(d-k) + d}|\det T_1|^{-1}
   \end{split}
\end{equation}
as before.  
Thus for $D_0 = M^{d\ell - \eta}$ we  obtain after dyadic decomposition
  \begin{equation}\label{cll}
  \mathcal{C}_{\ell,-}(M, D_0 ) \ll  \sup_{D \ll D_0}  D^{-1} M^{(\ell-1)(d-k) + d+\varepsilon} \mathcal{C}^{\dagger}_{ \ell}(M, D)
  \end{equation}
   where 
   \begin{align}
   \mathcal{C}^{\dagger}_{ \ell}(M, D)=\#\Big\{(x_{j,i})_{\ell \times k}:& \ x_{j,i}\in \mathbb{Z},\  x_{j,i}\ll M, \ |\det T_1|\asymp D,\\& T' \text{ is optimally arranged in the sense of \eqref{oa}}\Big\}.
   \end{align}

\begin{lemma}\label{lem3} Let $k \geq  d\geq 2$ and $k\geq \ell+1\geq 3$. We have
$$\mathcal{C}_{\ell,-}(M, D_0 ) \ll M^{k+\varepsilon}\Big( \frac{D_0}{M^{d\ell}} +   M^{ \ell(d + 1) - 1-k}\Big)$$ 
for any $\varepsilon>0$.
\end{lemma}

\textbf{Proof.} 
The quantity $\mathcal{C}_{\ell,-}(M, D_0 )$   does not change if we replace $T$ with $T'$ (since it depends only on $\det T_1$).  
For notational simplicity let us assume that $\ell$ is even, the case $\ell$ odd is analogous.  We introduce one more piece of notation: we write
$$\Delta_n = \begin{cases}  |\mathcal{T}^{\{1, \ldots, n\}}_{\{1, \ldots, n\}} |, & n  \text{ even} \\ \max\Big(|\mathcal{T}^{\{1, \ldots, n\}}_{\{1, \ldots, n\}}|, |\mathcal{T}^{\{1, \ldots, n\}}_{\{1, \ldots, n-1, n+1\}}|\Big), & n \text{ odd}.  \end{cases}$$

We now have 
  $$\mathcal{C}_{\ell, -}(M, D_0 ) \ll   \sup_{\substack{D_1, \ldots, D_{\ell}\\ D_{\ell} \ll D_0}}  D_{\ell}^{-1} M^{(\ell-1)(d-k) + d+\varepsilon} \mathcal{C}^\dagger_{\ell}(M, D_1, \ldots, D_{\ell})$$
where
\begin{align}
\mathcal C_\ell^\dagger(M, D_1, \ldots, D_\ell )=\#\Big\{(x_{j,i})_{\ell \times k}:&\ x_{j,i}\in \mathbb{Z}, \ \  x_{j,i}\asymp M,\  T' \text{ is optimally arranged},\\& \Delta_n\asymp D_n \text{ for } 1\leq n\leq \ell \Big\}.
\end{align}
   Since $T'$ is optimally arranged and $\rank T=\ell$, we have $\Delta_n>0$ for all $n\geq 1$.
  
We distinguish two cases when estimating $\mathcal C^\dagger_{\ell}(M, D_1,\dots, D_\ell)$:

1) Suppose $D_n \geq D_{n-1} M^{d-1}$ for  all $1 \leq n \leq \ell$. We choose the entries of $T'$ as follows.  We choose $x_{1,1}, x_{2,1}$ in $O(M^2)$ ways. Next we choose $x_{1,2}, x_{2,2}$ in $O( M^{2-d} D_2/D_1)$ ways using  Lemma \ref{lemma-b}a  with $k=2$ and $H = D_2$. In the same way we choose $x_{1,3}, x_{2,3}$ keeping in mind the condition \eqref{oa} with $n=2$.  Next we choose $x_{3,1}, x_{3,2}, x_{3,3}$ in $O( M^{3-d} D_3/D_2)$ ways using Lemma \ref{lemma-b}a with $k=3$. In the same way we choose $x_{4,1}, x_{4,2}, x_{4,3}$. Next we see that we can choose $x_{1,4}, \ldots, x_{4,4}$ in $O(M^{4-d}D_4/D_3)$ ways using Lemma \ref{lemma-b}a together with the fact that $$D_3\leq \max_{\substack{|\mathfrak c|=|\mathfrak r|=3\\\mathfrak c \subset\{1,2,3,4\}\\ \mathfrak r\subset \{1,2,3\}}}|\mathcal T_{\mathfrak r}^{\mathfrak c}|$$
 and the same bound applies to $x_{1,5}, \ldots, x_{4,5}$ (also using \eqref{oa} with $n=4$). We continue until we have chosen  $x_{1, \ell}, \ldots, x_{\ell, \ell}$ and $x_{1, \ell + 1}, \ldots, x_{\ell, \ell+1}$. The remaining $k- \ell - 1$ columns can be chosen randomly. Thus the total count for $\mathcal{C}^{\dagger}_{\ell}(M, D_1, \ldots, D_{\ell})$ is
$$\ll M^2    \Big(M^{2-d} \frac{D_2}{M^d}\Big)^2 \Big(M^{3-d} \frac{D_3}{D_2}\Big)^2 \cdots \Big(M^{\ell-d} \frac{D_{\ell}}{D_{\ell-1}}\Big)^2 M^{\ell(k-\ell - 1)} = D_{\ell}^2 M^{k\ell - 2 d \ell},$$
 and hence
\begin{equation}\label{Dnlarge}
\max_{\substack{D_{\ell} \ll D_0\\ D_n\geq D_{n-1}M^{d-1} \\\text{ for all }1\leq n\leq \ell}} D_\ell^{-1}M^{(\ell-1)(d-k)+d+\varepsilon}C^{\dagger}_{\ell}(M, D_1, \dots, D_\ell)
\ll D_0 M^{k - d\ell+\varepsilon}.
\end{equation}

2) Now we assume that $D_n \leq D_{n-1} M^{d-1}$ for some $1 \leq n \leq \ell$. Here we choose the upper left $n \times (n-1)$-block randomly. Next we apply Lemma \ref{lemma-b}a with $k = n$ to pick the first $n$ elements in the $m$-th column $x_{1, m}^d, \ldots, x_{n, m}^d$ for $n \leq m \leq k$ in $$\ll M^{n-1}\Big( \frac{D_n}{M^{d-1} D_{n-1}} + 1\Big) \ll M^{n-1}$$ ways, keeping in mind \eqref{oa}.  Finally we choose the $\ell - n$ bottom rows randomly.  Thus the total count for $\mathcal{C}^{\dagger}_{\ell}(M, D_1, \ldots, D_{\ell})$ can be bounded by
$$\ll M^{n(n-1)} \cdot M^{(k - n+1)(n-1)} \cdot M^{k(\ell - n)} = M^{k(\ell-1) + n-1} \leq M^{(k+1)(\ell-1)} $$
and hence
\begin{equation}\label{Dnsmall}
\max_{\substack{D_{\ell} \ll D_0\\ D_n\leq D_{n-1}M^{d-1} \\\text{ for some }1\leq n\leq \ell}} D_\ell^{-1}M^{(\ell-1)(d-k)+d+\varepsilon}C^{\dagger}_{\ell}(M, D_1, \dots, D_\ell) \ll  M^{\ell(d + 1) - 1+\varepsilon}.
\end{equation}
Combining \eqref{Dnlarge} and \eqref{Dnsmall}, we obtain the lemma. 

\subsection{Completion of the proof} In order to prove Theorem \ref{thm3}, we apply Proposition \ref{prop2} using Lemmas \ref{lem2} and \ref{lem3}. The continuous contribution $M^k \mathcal{G}_-(\eta)$ can be bounded as in Lemma \ref{lem3}, but the analysis is much easier (cf.\ e.g.\ \cite[Section 4]{ABR} or \cite[Appendix]{VdK2}).  The contribution $\det T_1 = 0$ is a Lebesgue null set, and so we are left with bounding the volume of (real) matrices satisfying the conditions in the definition of $\mathcal{C}^{\dagger}_{\ell}(M, D_0)$. This can be done as in Lemma \ref{lem3}, except that we  can drop the $+1$-term in Lemma \ref{lemma-b}a when estimating integrals as opposed to sums. This principle applies in general, and in the following sections we will not discuss the continuous contribution in Propositions \ref{prop2} and \ref{prop3}. 

\section{Proof of Theorem \ref{thm1}}

The proof is very similar to that of Theorem \ref{thm3}. 
To begin with, let $\ell, d \geq 2$, $k \geq \max(2\ell, d)$.   We remark  that $\mathcal{C}^{(2)}(M; k, \ell, d)$ as defined in \eqref{quantity2} is very similar to $\mathcal{C}(M; k, 2\ell, d)$   as defined in \eqref{quantity}, but with two important differences: the asterisk means something slightly different, and we have only $2\ell - 2$ conditions $q _{\bm\alpha}(\bm x_j) - q _{\bm\alpha}(\bm x_1), q _{\bm\alpha}(\bm y_j) - q _{\alpha}(\bm y_1) \ll M^{d-k}$, rather than $2\ell - 1$ conditions as in $\mathcal{C}(M;k,2\ell, d)$. 

We apply Proposition \ref{prop3}, and from the proof of Theorem \ref{thm3} we know already \eqref{prop2-bound}. We bound $\mathcal{C}^{(2)}_r(M)$, $2 \leq r \leq 2\ell-1$,  and $\mathcal{C}^{(2)}_{2\ell,-}(M, D_0)$ exactly as in Lemmas \ref{lem2} and \ref{lem3}.  The only difference is that for $ \text{rank}({\tt T}_1) = r$ we have  
 \begin{equation}\label{intnew}
  \mathcal{I}^{(2)}(\bm{x}, \bm{y}, M) \ll M^{(r - 2)(d-k) + 2d} |\det \tt T_1|^{-1}
\end{equation} 
  since  two of the $a_j$'s can be size $M^d$. Compared to \eqref{inte} this loses a factor of $M^k$ (which is admissible since our target in Proposition \ref{prop3} differs by a factor $M^k$ from Proposition \ref{prop2}).  It only remains to take care of the new meaning of the asterisk in \eqref{quantity2}. In \eqref{quantity} we could assume up front that $\bm{x}_1, \ldots, \bm{x}_k$ are pairwise linearly independent.    Suppose now that  some $\bm x_j$ equal some   multiples of some  $\bm y_i$. In this case we lose a factor $M^{k-1}$ and we win a factor $M^k$, as we explain now. The loss comes from the fact that in the linear combination \eqref{linear} we may not have the condition \eqref{two}, so that we cannot apply Lemma \ref{lemma-a}a. Instead, we apply the trivial bound to the equation \eqref{linear} and replace \eqref{r-1} with $O(M^r)$. Since we apply this for $k-1$ columns, we lose a factor $M^{k-1}$. On the other hand, if $\bm{y}_1$ (without loss of generality, after relabeling the indices and accordingly the definition of $\tt T$) and $\bm{x}_i$ are linearly dependent for some $i$, then in the decomposition 
 \begin{equation}\label{tdecomp} 
  {\tt T} = \left( \begin{matrix} {\tt T}_1 & {\tt T}_2\\ {\tt T}_3& {\tt T}_4 \end{matrix}\right) \in \Bbb{R}^{(r + (2\ell-r)) \times (r + (k-r))} 
  \end{equation}
corresponding to \eqref{decomp} , the matrix ${\tt T_1}$ does not contain the row $(y_{1,1}^d, \ldots, y_{1,k}^d)$, so we in fact have the better bound 
$$\mathcal{I}^{(2)}(\bm{x}, \bm{y}, M) \ll M^{(r - 1)(d-k) + d} |\det \tt T_1|^{-1}$$ which saves a factor $M^k$ compared to \eqref{intnew}. The omission of the row $(y_{1,1}^d, \ldots, y_{1,k}^d)$ does not cause any problems when we pass to ${\tt T}' = \left(\begin{matrix} T'(\bm x)\\ T'(\bm y)\end{matrix}\right)  \in \Bbb{R}^{2\ell \times k}$ as in \eqref{T'}, since $\bm{y}_1$ is linearly dependent on some $\bm{x}_i$, so by row operations we can still transform a row $(y_{j,1}^d - y_{1,1}^d, \ldots, y_{j,k}^d - y_{1,k}^d)$  to pure $d$-th powers $(y_{j,1}^d , \ldots, y_{j,k}^d)$. 
The rest of the argument is identical.

\section{The case $\ell=3$}
In this section we prove Theorem \ref{thm4} and one half of Theorem \ref{thm5} by improving Lemmas \ref{lem2} and \ref{lem3} in the special case $\ell = 3$ which we assume throughout this section. We recall the definitions \eqref{TT'} and \eqref{T'} of the matrices $T, T'$. The proof of Theorem \ref{thm4} will follow from Lemmas \ref{lem4} and \ref{lem5} together with Propositions \ref{prop1} and \ref{prop2}. The proof for \eqref{weak} in Theorem \ref{thm5}  follows from Lemmas \ref{lem6s} and \ref{lem6} together with Propositions \ref{prop1} and \ref{prop2}. 

\subsection{Small ranks} We start by considering the case  $\rank(T) < 3$ in which case necessarily we have $\rank(T) = 2$ as observed in the beginning of the proof of Theorem \ref{thm3}.  

\begin{lemma}\label{lem4} Let $d\geq 2$, $k \geq \max(d, 3)$. We have 
$\mathcal{C}_{2}(M; k, 3, d) \ll   M^{\varepsilon}(M^2 + M^{1 + \frac{2}{d} + 2d-k})$ for every $\varepsilon > 0$. 
\end{lemma}

\textbf{Proof.} Since $\rank(T')=\rank T= 2$,  we can write the last row   of $T'$ as a linear combination of the first two rows, i.e.\  
\begin{equation}\label{linear-new}
x^d_{3,i} =  \rho_{1}  x_{1, i}^d + \rho_{2} x_{2, i}^d,  \quad 1 \leq i \leq k,
\end{equation}
and by exchanging rows if necessary we may assume that $0 < |\rho_1|, |\rho_2| \leq 1$.  This is the analogue of \eqref{linear}. 
Since  $M^d\asymp q_{\bm \alpha}(\bm x_3)=\rho_1 q_{\bm \alpha}(\bm x_1)+\rho_2 q_{\bm \alpha}(\bm x_2)$ for any ${\bm \alpha}$ in the support of $F$, we conclude 
\begin{equation}\label{sum}
\rho_1 +\rho_2 =1+O(M^{-k}), 
\end{equation}
which implies 
\begin{align}
x_{3,i}^d-x_{1,i}^d=\rho_2(x_{2,i}^d-x_{1,i}^d)+O(M^{d-k}).
\end{align}
 We call a column \emph{typical} if all three elements are distinct. 
And we observe that a column of $T'$ is either typical or has identical entries. Indeed, if two of them are equal, say $x_{2, i} = x_{1, i}$, then \eqref{linear-new} and \eqref{sum} imply $x_{3, i} = (1 + O(M^{-k})) x_{2, i}$, and since $x_{2, i}, x_{3, i} \in \Bbb{N}$ we must have equality. If $x_{3,i}=x_{2,i}$ for some $i$, then $|\rho_1|\ll M^{1-k}$ and $|x_{3,j}^d-x_{2,j}^d|\ll M^{d+1-k}$ for all $j$, which is impossible if $x_{2,j}\not=x_{3,j}$. Thus we must have $x_{2,i}=x_{3,i}$ for all $i$, which is impossible.

We argue now that we can always arrange that at least two columns of $T'$, say the first two, must be typical and in addition the upper left $2\times 2$-block $T^{\natural}$ of $T'$ must be non-singular. The first claim is clear: if only the first column was typical, then $q_{\bm \alpha}(\bm{x}_1) - q_{\bm \alpha}(\bm{x}_2) = \alpha_1(x_{1,1}^d - x_{2,1}^d) \gg M^{d-1}$ for $\alpha_1  \gg 1$, a contradiction. And if the first two components of all typical columns were linearly dependent, then  $x_{1,i} = x_{2,i}$ when the  $i$-th column is  non-typical and $x_{1,i} = \beta x_{2,i}$ with $\beta \not= 1$ when the $i$-th column is typical, and again $$q_{\bm \alpha}(\bm{x}_1) - q_{\bm \alpha}(\bm{x}_2) = \sum_{i \text{ typical}}\alpha_i(x_{1,i}^d - x_{2,i}^d) \gg M^{d-1}$$ for $\alpha_i \gg 1$ by the same argument, a contradiction. 

Let us in addition order the first two columns such that   
\begin{align}\label{Xdef}
Y\asymp  |x_{2, 2} - x_{1, 2}| \leq |x_{2, 1} - x_{1, 1}| \asymp X
\end{align}  for some parameter $1 \leq Y\leq X \ll M$. Let $g=(x_{3,2},x_{2,2},x_{1,2})$ and fix it.  For $i = 1, 2$ we can define
\begin{equation}\label{qidef}
q_{i}=\frac{x_{3,i}^d-x_{1, i}^d}{x_{2, i}^d-x_{1,i}^d} = \rho_2 + O\Big(\frac{M^{1-k}}{|x_{2, i} - x_{1, i}|}\Big) = q_1 +  O\Big(\frac{M^{1-k}}{X}\Big)=q_2+O\Big(\frac{M^{1-k}}{Y}\Big).
\end{equation} 
We choose the first column in $O(M^2X)$ ways. This determines $q_1$, which in turn   determines  $q_2$  as a reduced fraction with denominator of size $O(YM^{d-1}/g^d)$ in at most $\ll 1+ M^{1-k}Y^{-1} (YM^{d-1})^2 g^{-2d} \ll 1+ M^{2d-k}g^{-2d}$ ways, using $Y\ll X$. Since the second column is typical, we have $q_2 \not= 0, 1$. Then we can apply Lemma \ref{lemma-a}b to the first equation in \eqref{qidef} with $i=2$ to determine the second column in $\ll (M/g)^{2/d+\varepsilon}$ ways. We conclude that the first two columns can be chosen in $$\ll M^{\varepsilon} \sum_{g\ll M}M^2 X \Big(\frac{M}{g}\Big)^{2/d} \Big(1 + \frac{M^{2d-k}}{g^{2d}}\Big)\ll M^{2+2/d+\varepsilon}X(M^{(d-2)/d}+M^{2d-k})$$ ways. 

At this point we can compute $\rho_1, \rho_2$ (both different from $0, 1$) from the non-singular matrix equation $$(T^{\natural})^{\top}\left(\begin{matrix}\rho_1\\ \rho_2\end{matrix}\right) = \left(\begin{matrix}x_{3,1}^d\\ x_{3,2}^d\end{matrix}\right).$$ Hence for each remaining column, we can apply Lemma \ref{lemma-a}a to \eqref{linear-new} to get at most $O(M^{1+\varepsilon})$ solutions for each remaining column. While still restricting our count to elements satisfying $|x_{2, 1} - x_{1, 1}| \asymp X$, we obtain at most 
\begin{align}\label{c2M}
M^{\varepsilon} X (M^3   + M^{2 + \frac{2}{d} + 2d-k}) M^{k-2}.
\end{align} 
choices for all columns. 
Before we substitute this into \eqref{clr}, we also improve \eqref{inte} a bit. Namely, by Cramer's rule applied to the function $F({\bm \alpha}) = F(T_1^{-1}(\bm a - T_2 \bm \alpha_2), {\bm \alpha}_2)$ we may bound the integral over $a_1$ by $\ll |\det T_1|/|x_{2,1}^d - x_{1,1}^d| \ll |\det T_1|/(M^{d-1}X)$ rather than $M^d$, and hence
$\mathcal{I}(\bm{x}, M) \ll M^{d-k}/ (M^{d-1}X)$. 
 Putting this together with \eqref{c2M}, we obtain
$$\mathcal{C}_{2}(M) \ll  M^{\varepsilon} \sup_{1 \leq X \ll M}  \frac{M^{d-k}}{M^{d-1} X} X (M^3   + M^{2 + \frac{2}{d} + 2d-k}) M^{k-2}, $$
and hence the lemma.

\subsection{Full rank} Next we consider the situation  when $\text{rank}(T) = 3$ and tighten the argument of Lemma \ref{lem3}. We
  assume that $T$ is ``optimally arranged'' in the following way. Among all differences $|x_{j,i} - x_{j' ,i}|$, $1 \leq i \leq k$, $j \not=  j' \in \{1, 2, 3\}$, we choose the largest and assume without loss of generality that $i = 1$ (by re-ordering columns). Next we relabel the rows if necessary and assume that 
 \begin{equation}\label{sizeX}
 |x_{2,1}-x_{3,1}|\leq |x_{1,1}-x_{3,1}|\leq |x_{1,1}-x_{2,1}|.
 \end{equation}
Now we choose the second and third column appropriately to ensure 
  \begin{equation}\label{oa1}
  \begin{split}
   & |x_{2,1} - x_{1,1}| \geq |x_{2, i} - x_{1, i}|, \quad i \geq 1,\\
   & |\tilde{\mathcal{T}}^{\{1, 2\}}_{\{2, 3\}} | \geq |\tilde{\mathcal{T}}^{\{1, i\}}_{\{2, 3\}} |,  \quad   i \geq 2,\\
   & |\tilde{\mathcal{T}}^{\{1, 2, 3\}}_{\{1, 2, 3\}} | \geq |\tilde{\mathcal{T}}^{\{1, 2, i\}}_{\{1, 2, 3\}} |, \quad  i \geq 3,
    \end{split}
  \end{equation}
 where  $\tilde{\mathcal T}_{\mathfrak r}^{\mathfrak c}$ is the determinant of the submatrix of $T$ consisting of the columns in $\mathfrak c$ and rows in $\mathfrak r$ (the first condition is automatic for our choice of the first  column and \eqref{sizeX}).

We  decompose  $T=(T_1 \ T_2)$, where $$T_1 = \left(\begin{matrix} x_{1,1}^d & x_{1,2}^d & x_{1,3}^d\\ x_{2,1}^d - x_{1,1}^d & x_{2,2}^d - x_{1,2}^d & x_{2,3}^d - x_{1,3}^d\\x_{3,1}^d - x_{1,1}^d & x_{3,2}^d - x_{1,2}^d & x_{3,3}^d - x_{1,3}^d \end{matrix}\right)$$ is   regular   by the third condition in \eqref{oa1} and $\rank T=3$. We have $x_{2,1}\not=x_{1,1}$ by the first condition in \eqref{oa1} and $\rank T=3$. We also have $|\tilde{\mathcal T}_{\{2,3\}}^{\{1,2\}}|\not=0$ by the second condition in \eqref{oa1} and $\rank T=3$.
  As in the previous proof, we can can use Cramer's rule to slightly improve the bound for $\mathcal{I}(\bm{x}, M)$ in \eqref{quantity-new}  to 
  \begin{equation}\label{sizeint}
    \mathcal{I}(\bm{x}, M) \ll   M^{2(d-k)} |\tilde{\mathcal T}_{\{2,3\}}^{\{1,2\}}|^{-1},
 \end{equation}    
       and correspondingly we have
      $$ \mathcal{C}_{3,-}(M, D_0) = \sup_{\substack{D \ll D_0\\ Z \leq X \ll M\\\Delta \ll M^{2d}\\ X, \Delta>0}} \frac{M^{2(d-k)}}{\Delta} \mathcal{C}^{\dagger}_{3}(M, D, \Delta; X, Z)$$
      where 
      \begin{align}
      \mathcal{C}^{\dagger}_{3}(M, D, \Delta; X,Z)
      =\#\Big\{(x_{j,i})_{3\times k}:\ & x_{j,i}\in \mathbb{Z}, \  x_{j,i}\asymp M, \  |\det T_1| \ll D, \ |\tilde{\mathcal T}_{\{2,3\}}^{\{1,2\}}| \asymp \Delta \not= 0,\\&|x_{2,1} - x_{1,1}| \asymp X \not= 0,\  |x_{2,1} - x_{3,1}| \asymp Z\Big\}.
   \end{align}

      \begin{lemma}\label{lem5} For $d\geq 2, k \geq 4$ we have 
      	$$\mathcal{C}_{3,-}(M, D_0) \ll M^{k+\varepsilon} \Big( \Big(  \frac{D_0 }{M^{3d}}\Big)^{1/5} + M^{2d +\frac{2}{d} -\frac{3}{2}k} + M^{-2/3}\Big)$$
      	      for every $\varepsilon>0$.
      \end{lemma}

      \textbf{Proof.} We estimate $\mathcal{C}^{\dagger}_{3}(M, D, \Delta; X, Z)$ in various ways. 
      
      A first count is very simple: We pick the first  two columns randomly in $O(M^6)$ ways. For the $i$-column we obtain an inequality from the third condition in \eqref{oa1} which reads
      $$|A x_{1,i}^d + B (x_{2,i}^d - x_{1,i}^d) + C(x_{3,i}^d - x_{1,i}^d) |\ll D$$
      for certain numbers $A, B, C$ with $A=|\tilde{\mathcal T}_{\{2,3\}}^{\{1,2\}}|\asymp \Delta$. Since $\max(|A - B - C|, |B|, |C|) \gg |A| \asymp \Delta>0$, we can apply 
      Lemma \ref{lemma-b}a with $k=3$ to see that there are at most  $\ll M^{3-d} D/\Delta + M^2$ choices for the $i$-th column. Alternatively, we can bound this number by $O(M^3)$ trivially. Thus we obtain 
      \begin{equation}\label{bound1}
      \begin{split}
      \frac{M^{2(d-k)}}{\Delta}  \mathcal{C}^{\dagger}_{3}(M, D, \Delta; X, Z)& \ll  \frac{M^{2(d-k)}}{\Delta} M^{6} \Big( \min\Big( \frac{M^{3-d}D}{\Delta} , M^3\Big)+ M^2\Big)^{k-2}\\
      & \ll \frac{M^{2(d-k)}}{\Delta} M^{6} \Big(\frac{M^{3-d}D}{\Delta}  M^{3(k-3) }+   M^{2(k-2)}\Big)\\
      &\ll   \frac{M^{d+k}D }{\Delta^2} + \frac{M^{2 + 2d}}{\Delta} .
      \end{split}
      \end{equation}
      This is useful if $\Delta$ is big.     
      
      When $\Delta$ is smaller we can study the expression of $\tilde{\mathcal T}_{\{2,3\}}^{\{1,2\}}$ using Lemma \ref{lemma-b}c. We have 
      \begin{equation}\label{d2}
      |\tilde{\mathcal T}_{\{2,3\}}^{\{1,2\}}| =  |(x_{3,1}^d-x_{2,1}^d)x_{1,2}^d+(x_{1,1}^d-x_{3,1}^d)x_{2,2}^d+(x_{2,1}^d-x_{1,1}^d)x_{3,2}^d|. 
      \end{equation}
     We first dispense the case when $Z=0$. In this case, we see that \eqref{d2} becomes 
      \begin{align}
      |(x_{1,1}^d-x_{3,1}^d)(x_{2,2}^d-x_{3,2}^d)|\asymp \Delta
      \end{align}
      Since $X>0$, we must have $0<|x_{1,1}- x_{3,1}|\asymp Y\ll X$. 
      Since $\Delta>0$, we see that the number of choices for the second column can be bounded by $O\Big( M^2  \Delta/(M^{2d-2}Y)\Big)$. Using $O(MY)$ for the number of choices for the first column and the trivial bound $O(M^3)$ for all other columns we have 
      \begin{equation}\label{Z=0}
      \frac{M^{2(d-k)}}{\Delta}\mathcal C_3^\dagger (M,D, \Delta;X,0)\ll \frac{M^{2(d-k)}}{\Delta} \sup_{0<Y\ll M}M Y \frac{M^2\Delta}{M^{2d-2}Y} M^{3(k-2)}\ll M^{k-1},
      \end{equation}
      which is admissible.
      
       From now we  assume $Z>0$.  By \eqref{sizeX}   we have $O(MXZ)$ choices for the first column, and from \eqref{d2} we see that   the number  of $x_{1,2}, x_{2,2},x_{3,2}$ such that $|\tilde{\mathcal T}_{\{2,3\}}^{\{1,2\}}|\ll \Delta$ is, by Lemma \ref{lemma-b}c, bounded by 
      \begin{align}
      &M^{\varepsilon}\Big( \frac{M^{3-d}\Delta}{(XZ)^{1/2}M^{d-1}}+ \frac{M^{3-d/2}\Delta^{1/2}}{(XM^{d-1})^{1/2}}+\frac{M^{5/2-d/2}\Delta^{1/2}}{(ZM^{d-1})^{1/2}}+M^{3/2}\Big)\\
      &= M^\varepsilon \Big( \frac{M^{4-2d}\Delta}{(XZ)^{1/2}}+\frac{M^{7/2-d}\Delta^{1/2}}{X^{1/2}}+ \frac{M^{3-d}\Delta^{1/2}}{Z^{1/2}}+M^{3/2}\Big).
      \end{align}
      By the second condition in \eqref{oa1}, we can apply this reasoning for every $i$-th column with $i\geq 2$. Alternatively, we can choose the $i$-th column trivially in $O(M^3)$ ways. In this way we obtain for $k\geq 4$ that
      \begin{align}
      &\frac{M^{2(d-k)}}{\Delta}  \mathcal{C}^{\dagger}_{3}(M, D, \Delta; X, Z) 
      \\& \ll  \frac{M^{2(d-k)+\varepsilon}}{\Delta}  M XZ \Big( \min\Big(\frac{M^{4-2d}\Delta}{(XZ)^{1/2}}+\frac{M^{\frac{7}{2}-d}\Delta^{1/2}}{X^{1/2}}+ \frac{M^{3-d}\Delta^{1/2}}{Z^{1/2}}, M^3\Big) + M^{\frac{3}{2}}\Big)^{k-1}\\
      & \ll \frac{M^{2(d-k)+\varepsilon}}{\Delta}  M XZ \Big( \Big(\frac{ M^{4-2d}  \Delta}{  (XZ)^{1/2}}\Big)^2 M^{3(k-3)} + \Big(\frac{M^{7/2-d}\Delta^{1/2}}{X^{1/2}}\Big)^3M^{3(k-4)}\\  \label{bound2}
      & \quad\quad\quad\quad\quad\quad\quad\quad\quad\quad+ \frac{M^{6-2d}\Delta}{Z}M^{3(k-3)}+M^{3(k-1)/2}\Big)\\ 
      & \ll M^{k+\varepsilon}\Big( \frac{\Delta}{M^{2d}}+ \frac{\Delta^{1/2}}{M^d}+\frac{1}{M}+ \frac{M^{(3 + 4 d - 3 k)/2}}{\Delta}\Big).
      \end{align}
      This would suffice for $k> \frac{1}{3}(4d+3)$ and we improve this by a third argument under the general assumption 
      \begin{equation}\label{D2bound}
      0 < \Delta 
      \ll M^{2d-2} 
      \end{equation}
       with a sufficiently small implied constant. 
       We write the second condition in \eqref{oa1} for the determinant $|\tilde{\mathcal{T}}^{\{1, i\}}_{\{2, 3\}} |$ with $i\geq 2$ as
      \begin{align}\label{D2inequal2}
      |(x_{2,1}^d-x_{1,1}^d)(x_{3,i}^d-x_{1,i}^d)-(x_{3,1}^d-x_{1,1}^d)(x_{2,i}^d-x_{1,i}^d)| =:   D_i  \ll D_2 \asymp \Delta.
      \end{align}
       
      As before we call the $i$-th column typical if the set $C_i = \{x_{1, i}, x_{2, i}, x_{3, i}\}$ has cardinality 3. We first observe that as long as $D_i \not= 0$, the case $\#C_i < 3$ or $\#C_1 < 3$ cannot occur, since otherwise \eqref{D2bound} would be violated. This applies in particular for $i = 2$ and so the first two columns are typical. 
      On the other hand, if $D_i = 0$ for some $i > 2$, then either the $i$-th column is typical or we must have  $x_{2,i} = x_{1,i} = x_{3,i}$ since we know that the first column is typical and $D_i\ll \Delta$ satisfies \eqref{D2bound}, or we have $$x_{3,i}^d-x_{1,i}^d=q_1(x_{2,i}^d-x_{1,i}^d)$$ with $q_1\not=\{0,1\}$ and Lemma \ref{lemma-a}(a) gives $O(M^{1+\varepsilon})$ choices for the $i$-th column. 
      
      We first choose the first column in $\ll MX^2 \ll M^3$ ways. Let us now consider the $i$-th column for some $i \geq 2$. If this column is not typical, then we can choose it in $O(M^{1+\varepsilon})$ ways. Suppose now that it is typical, then we can assume 
      $|x_{2,i}-x_{1,i}|\asymp Y $ with $0<Y \leq X$ by the first condition in \eqref{oa1}. Now we can re-write \eqref{D2inequal2} as 
      \begin{equation}\label{q1qi}
      |q_1 - q_i| 
      \ll  \frac{\Delta}{M^{2(d-1)}X Y }, \quad q_i =  \frac{x_{3,i}^d-x_{1,i}^d}{x_{2,i}^d-x_{1,i}^d} \not\in \{ 0, 1\}. 
      \end{equation}
      Let $g_i=(x_{1,i},x_{2,i},x_{3,i})$ and fix it. Since $q_1$ is already fixed by the first column, the number of choices for $q_i$ as a fraction with denominator $\asymp YM^{d-1}g_i^{-d}$ is at most  
      \begin{equation}\label{qichoice}
      1+ \frac{\Delta}{M^{2(d-1)}XY}(YM^{d-1}g_i^{-d})^2 \ll 1 + \Delta \frac{Y}{Xg_i^{2d}} \ll 1+\frac{\Delta}{g_i^{2d}}. 
      \end{equation}
      Once $g_i, q_i$ are fixed with $q_i\not\in\{0,1\}$, we can apply Lemma \ref{lemma-a}b to \eqref{q1qi} to see that the $i$-th column can be chosen in  $O((M/g_i)^{2/d+\varepsilon})$ ways. Summing over $g_i\ll M$ we see that the number of choices for the $i$-th column can be bounded by 
      \begin{align}
      M+\sum_{g_i\ll M}\Big(1+\frac{\Delta}{g_i^{2d}}\Big)\Big(\frac{M}{g_i}\Big)^{2/d+\varepsilon}\ll M^{1+\varepsilon}+\Delta M^{2/d+\varepsilon}. 
      \end{align}
      Thus we obtain altogether
      \begin{align}\label{bound3}
      \frac{M^{2(d-k)}}{\Delta}  \mathcal{C}^{\dagger}_{3}(M, D, \Delta; X, Z)  &\ll  \frac{M^{2(d-k)+\varepsilon}}{\Delta}  M^3(M+\Delta M^{2/d})^{k-1}  \\
      &\ll  M^{\varepsilon}( M^{2d - k + 2  } + \Delta^{k-2} M^{2d + 3 - 2k + 2(k-1)/d})
      \end{align}
      under the assumption \eqref{D2bound}. 
      
We now combine \eqref{Z=0}, \eqref{bound2} and \eqref{bound3} as follows: if $\Delta \gg M^{3/2 - 2/d}$ we apply \eqref{Z=0} and \eqref{bound2}, otherwise we can use \eqref{bound3} since \eqref{D2bound} is satisfied. In this way we obtain
$$ \frac{M^{2(d-k)}}{\Delta}  \mathcal{C}^{\dagger}_{3}(M, D, \Delta, X, Z)  \ll M^{k+\varepsilon} \Big( M^{2(d+1-k)} + M^{2d - 3k/2 + 2/d } + \frac{1}{M}+ \frac{\Delta^{1/2}}{M^{d}}+\frac{\Delta}{M^{2d}}\Big).$$
We can drop the first term for $k \geq 4$ and combine the rest with \eqref{bound1}. Using appropriate geometric means, we have
\begin{displaymath}
\begin{split}
\min\Big( \frac{M^d D}{\Delta^2} + \frac{M^{2+2d-k}}{\Delta}, \frac{\Delta}{M^{2d}}+ \frac{\Delta^{1/2}}{M^d}\Big)&  \ll \frac{D^{1/3}}{M^d} + M^{1 - k/2} + \frac{D^{1/5}}{M^{3d/5}} + M^{(2-k)/3}\\
& \ll \frac{D^{1/5}}{M^{3d/5}} + M^{(2-k)/3} \leq \frac{D^{1/5}}{M^{3d/5}} + M^{-2/3}
\end{split}
\end{displaymath}
for $k \geq 4$ and $D \ll M^{3d}$. 
In this way we obtain
$$ \frac{M^{2(d-k)}}{\Delta}  \mathcal{C}^{\dagger}_{3}(M, D, \Delta, X, Z)  \ll M^{k+\varepsilon} \Big(  M^{ 2 d -3k/2+2/d} + M^{-2/3} + \frac{D^{1/5}}{M^{3d/5}}\Big).$$
This completes the proof.\\

Lemmas \ref{lem4} and \ref{lem5} suffice as soon as $k > \frac{4d}{3}+\frac{4}{3d}$ in which case we need $k \geq 4$ when $d=2$. In the following we improve Lemma \ref{lem4} and Lemma \ref{lem5} to allow $k= 3$ when $d=2$.

\subsection{Small ranks: the case $d=2$}  
We start with a variation of Lemma \ref{lem4} for $d=2$. 
  
  \begin{lemma}\label{lem6s} We have 
  	$\mathcal{C}_{2}(M; 3, 3, 2) \ll   M^{7/3+\varepsilon}$ for any $\varepsilon>0$.
  \end{lemma}
  
  \textbf{Proof.} 
  We assume throughout this proof $(k, \ell, d) = (3, 3, 2)$. 
  As in the proof of Lemma \ref{lem4} we can assume that the first two columns are typical and define \eqref{qidef} for $i=1,2$. We write
  \begin{equation}\label{qab}
  q_i = \frac{a_i}{b_i}, \quad  (a_i, b_i) = 1, \quad b_i > 0,\quad q_i\not\in\{0,1\}, \quad i=1,2.
  \end{equation}
  We recall the definition of $X$ in \eqref{Xdef}, which gives  $O(M^2X)$ choices for the first column. Once the first column is chosen, we have $q_1$ is determined and then we decompose $b_1$ into dyadic ranges with $b_1 \asymp B \ll MX$.  
  
  Next we consider the second column. Let $g = (x_{1, 2}, x_{2, 2}, x_{3, 2})$ and fix it.  
  Let us first assume that $q_2 \not= q_1$. In this case we observe that $|q_2 - q_1| \geq 1/(b_1b_2)$, so \eqref{qidef} implies  $b_2 \gg M^2 Y/B$. We also have $b_2\ll MY/g^2$ from \eqref{Xdef} and thus we must have $g^2\ll B/M\ll X$ whenever $b_2$ exists. 
  We determine $q_2$ from \eqref{qidef} in $\ll 1+M^{-2}Y^{-1} (YM/g^2)^2 = 1+ Y/g^4$. Now \eqref{qidef} yields a non-degenerate  ternary quadratic form in $x_{j,2}/g, 1\leq j\leq 3$ with determinant $|a_2b_2(b_2 - a_2)| \gg |b_2|^2$. Since $(a_2, b_2) = 1$, Lemma \ref{lemma-a}c shows that the triple $x_{j,  2}/g$, $1 \leq j \leq 3$, can be chosen in at most 
  \begin{equation}\label{HBbound}
  \ll M^{\varepsilon}\Big(1 + \frac{M}{g|b_2|^{2/3}}\Big) \ll M^{\varepsilon}\Big(1 + \frac{M B^{2/3}}{g(M^2 Y)^{2/3}}\Big)
  \end{equation}
  ways. Thus the second column (with fixed $g$) can be chosen in 
  \begin{align}
&\ll M^{\varepsilon}\sum_{g\ll \sqrt{X}}\Big(1+\frac{Y^4}{g^4}\Big)\Big(1 + \frac{M B^{2/3}}{g(M^2 Y)^{2/3}}\Big)\\
&\ll M^\varepsilon\Big(\sqrt{X}+Y+M^{1/3}X^{2/3}Y^{-2/3}+M^{1/3}X^{2/3}Y^{1/3}\Big)\ll M^{4/3+\varepsilon}
  \end{align}
  ways if $q_2 \not= q_1$.  On the other hand, if $q_2 = q_1$, then simply by Lemma \ref{lemma-a}a we can choose the triple $x_{j, 2}/g$, $1 \leq j \leq 3$, in $O(M^{1+\varepsilon}/g)$ ways. 
  
  Finally, once the first and second columns are determined, as in the proof of Lemma \ref{lem4} the last column is determined in $O(M^{1+\varepsilon})$ ways. For $d=2, \ell=k=3$ we conclude that 
  \begin{displaymath}
  \begin{split}
  \mathcal{C}_{ 2}(M)&  \ll M^{\varepsilon} \sum_{g \ll M} \sup_{\substack{0< X \ll M\\ B \ll XM}}  \frac{1}{M^2X}  M^2X\Big( \frac{M}{g} + M^{4/3}\Big) M \ll M^{7/3 + \varepsilon}. 
  \end{split}
  \end{displaymath}
  This completes the proof. 
  
\subsection{Full rank: the case $d=2$}  We proceed with a variation of Lemma \ref{lem5} for $d=2, \ell=k=3$.

 \begin{lemma}\label{lem6} If $d= 2$, $\ell=k = 3$, then
$$ \mathcal{C}_{ 3,-}(M, D_0) \ll  M^{3+\varepsilon} \Big( \frac{D_0^{1/3}}{ M^2} + M^{- \frac{1}{10}}\Big) $$
for every $\varepsilon>0$.
\end{lemma}

\textbf{Proof.}  We want to improve \eqref{bound3} if $d=2$. Again we see that the assumption \eqref{D2bound} implies that the first two columns are typical and the third column is typical unless $D_3=0$ (with the notation as in \eqref{D2inequal2}), in which case $x_{1,3} = x_{2,3} = x_{3,3}$. Therefore we see that  \eqref{q1qi} is well defined for $i=1,2$ and we can write
$$q_i=\frac{a_i}{b_i}, \quad  (a_i, b_i)=1, \quad b_i > 0, \quad q_i\not\in\{0,1\}\quad i=1,2.$$
Recall $Y\asymp |x_{2,2}-x_{1,2}|\leq |x_{2,1}-x_{1,1}|\asymp X$.  Define $g=(x_{1,2}, x_{2,2},x_{3,2}) \ll M$ and fix it. From \eqref{D2inequal2} we have $g^2 \mid D_2\asymp \Delta$ and thus $g \leq \Delta^{1/2}$. 
 We fix the first column in 
 $O(MX^2)$ ways and this determines $q_1$.
 Next we consider the second column. Since $D_2\asymp \Delta>0$, we have $q_2\not=q_1$. Then from $|q_2 - q_1| \geq (b_1b_2)^{-1}$ and \eqref{q1qi} with $i = 2$, $d=2$ we must have 
\begin{align}
b_2 \geq \frac{M^2XY}{\Delta b_1}
\end{align}
whenever $q_2$ exists.
Since $q_2\not=0,1$, we see that $x_{1,2},x_{2,2},x_{3,2}$ satisfy a non-degenerate ternary quadratic equation with determinant $|a_2b_2(b_2-a_2)| \gg b_2^2$. Since $(a_2, b_2)=1$, Lemma \ref{lemma-a}c shows that we can choose the triple $x_{j,2}/g, 1\leq j\leq 3$,  in at most 
\begin{align}\label{x2q2}
\ll M^\varepsilon \Big(1+\frac{M}{g|b_2|^{2/3}}\Big)\ll M^\varepsilon \Big( 1+\frac{M|\Delta b_1|^{2/3}}{g(M^2 XY)^{2/3}}\Big) \ll M^{\varepsilon}\Big( 1+\frac{M^{1/3} \Delta^{2/3}  }{gY^{2/3}}\Big)
\end{align} 
ways for fixed $q_2\not=q_1$ and fixed $g$ (using $|b_1|\ll MX$). As in \eqref{qichoice}, the number of $q_2$ can be bounded by 
\begin{align}
1+\frac{\Delta }{M^{2} XY}\Big(\frac{YM}{g^2}\Big)^2\ll 1+\frac{\Delta Y}{g^4X} \ll \Delta.
\end{align} 
which together with \eqref{x2q2} gives a total contribution 
\begin{align}\label{x2bound1}
M^\varepsilon \sum_{g\ll \Delta^{1/2}}\Big(1+\frac{\Delta Y}{g^4X}\Big) \Big( 1+\frac{M^{1/3} \Delta^{2/3}  }{gY^{2/3}}\Big) \ll M^{\varepsilon}\Big( \Delta^{2/3}M^{1/3} + \Delta + \frac{M^{1/3} \Delta^{5/3}}{X^{2/3}}\Big)
\end{align}
for the number of choices for the second column.
We can use \eqref{x2bound1} for the number of choices for the third column when $D_3\not=0$ and $O(M)$ when $D_3=0$. Therefore we have
\begin{equation}\label{bound4}
\begin{split}
& \frac{1}{M^2\Delta}\mathcal C_3^\dagger(M, D, \Delta;X, Z)\\
&\ll \frac{M^\varepsilon}{M^2\Delta}   MX^2 \Big(\Delta^{2/3}M^{1/3}+\Delta+\frac{M^{1/3} \Delta^{5/3} }{X^{2/3}}\Big)\Big( M + \Delta + \frac{M^{1/3} \Delta^{5/3} }{X^{2/3}}\Big)\\
&\ll M^{\varepsilon}(\Delta ^{7/3} M^{1/3} + \Delta^{4/3} M + M^2 + \Delta^{2/3} M^{5/3} ) \ll M^{\varepsilon}(\Delta ^{7/3} M^{1/3} + M^{11/5} ).
\end{split}
\end{equation}

We also need to improve \eqref{bound2} for $d=2,\ell= k=3$, which was based on Lemma \ref{lemma-b}c.  Here we can use the special case $d=2$ of Lemma \ref{lemma-b}c in \eqref{d2} together with \eqref{sizeX} for the second and third columns and obtain 
\begin{equation}\label{k=3bound2}
\begin{split}
& \frac{1}{M^2\Delta} \mathcal C_3^\dagger (M, D, \Delta; X, Z)\\
&\ll \frac{M^{\varepsilon}}{M^2\Delta} \sup_{0<Z \leq Y \leq X}  M YZ \Big( \frac{M^2 \Delta}{YM}+\frac{M \Delta}{ZM}+M^2\Big)\Big(M+\frac{\Delta}{X M}\Big)\\
&\ll M^\varepsilon\Big( \frac{\Delta}{M}+ \frac{M^{4}}{\Delta}+M^2\Big).
\end{split}
\end{equation}
We now finish the argument exactly as before, but by replacing \eqref{bound3} with \eqref{bound4} and \eqref{bound2} by \eqref{k=3bound2}. Combining \eqref{bound4} and \eqref{k=3bound2} using suitable geometric means gives 
$$\frac{1}{M^2\Delta} \mathcal C_3^\dagger (M, D, \Delta, X, Z)\ll M^{\varepsilon}\big(\Delta M^{-1}+ M^{3-\frac{1}{10}} + M^{11/5}\big) \ll M^{\varepsilon}\big(\Delta M^{-1}+ M^{3 - \frac{1}{10}}\big).$$
Combining this with  \eqref{bound1} and \eqref{Z=0} completes the proof of the lemma. 

\section{The case $\ell = 2$}
In this final section we consider $L^2$-convergence for the pair correlation function and prove Theorem \ref{thm2} and the other half of Theorem \ref{thm5}.  Throughout this section we assume $\ell = 2$. 

 As in the proof of Theorem \ref{thm1}, we consider the matrix
$$\tt T = \left(\begin{matrix} x_{1,1}^d & \ldots & x_{1,k}^d\\ x_{2,1}^d - x_{1,1}^d & \ldots & x_{2,k}^d - x_{1,k}^d\\ y_{1,1}^d & \ldots & y_{1,k}^d\\ y_{2,1}^d - y_{1,1}^d & \ldots & y_{2,k}^d - y_{1,k}^d \end{matrix}\right).$$
For the proof of Theorem \ref{thm2} we assume $k \geq 4$. We can arrange the columns so that 
  \begin{equation}\label{oa2}
  \begin{split}
  &   |x_{2,1} - x_{1,1}| \geq |x_{2, i} - x_{1, i}|,\quad i\geq 1, \quad |x_{2,1} - x_{1,1}| \geq  |y_{2,1}-y_{1,1}|, \\
  & |\check{\mathcal {T}}^{\{1, 2\}}_{\{2, 4\}} | \geq |\check{\mathcal{T}}^{\{1, i\}}_{\{2, 4\}} |,  \quad   i \geq 2,\\
  & |\check{\mathcal {T}}^{\{1, 2,3\}}_{\{1,2, 4\}} | \geq |\check{\mathcal{T}}^{\{1,2,i\}}_{\{1,2, 4\}} |,  \quad   i \geq 3,\\
   & |\check{\mathcal {T}}^{\{1, 2,3,4\}}_{\{1,2, 3,4\}} | \geq |\check{\mathcal{T}}^{\{1,2,3,i\}}_{\{1,2,3, 4\}} |,  \quad   i \geq 4,
  \end{split}
  \end{equation}
  where  $\check{\mathcal T}_{\mathfrak r}^{\mathfrak c}$ is the determinant of the submatrix of $\tt T$ consisting of the columns in $\mathfrak c$ and rows in $\mathfrak r$.  We partition the count into 
\begin{equation}\label{XYDef}
 |x_{2,1} - x_{1,1}| \asymp X,\  |y_{2,1}-y_{1,1}|\asymp Y, \  |\check{\mathcal {T}}^{\{1, 2\}}_{\{2, 4\}} |  \asymp D_2, \  |\check{\mathcal {T}}^{\{1, 2,3\}}_{\{1,2, 4\}} |  \asymp D_3, \  |\check{\mathcal {T}}^{\{1, 2,3,4\}}_{\{1,2, 3,4\}} |  \asymp D_4.
\end{equation}
for parameters $X>0$ and $Y,D_2, D_3, D_4 \geq 0$. 

\subsection{Small ranks} We start with an upper bound for
$\mathcal{C}^{(2)}_{r}(M;k, \ell, d)$ as in \eqref{cr2def} for $r=2,3$. 
\begin{lemma}\label{c2rank2}
Let $ d\geq 2$, $k \geq \max(d, 4)$.	We have 
	$$ \mathcal C_{2}^{(2)}(M; k, 2, d)\ll M^{k+2}+M^{6+\varepsilon}$$
for any $\varepsilon>0$.
\end{lemma}
\textbf{Proof.}
The only way $\rank({\tt T}) = 2$ can happen is if $\bm{y}_1$ and $\bm{y}_2$ are linearly dependent on $\bm{x}_1$, $\bm{x}_2$. 
So we can write 
\begin{equation}\label{rhosigma}
y_{1,i}^d = \rho_1 x_{1,i}^d + \rho_2 x_{2,i}^d, \quad y_{2,i}^d = \sigma_{1} x_{1,i}^d + \sigma_{2} x_{2,i}^d.
\end{equation}
for some $\rho_i, \sigma_i\in \mathbb{R}$.
In the notation of \eqref{tdecomp} the change of variable matrix is the upper left $2$-by-$2$ block 
$${\tt T}_1 = \left(\begin{matrix} x_{1,1}^d &   x_{1,2}^d\\ x_{2,1}^d - x_{1,1}^d & x_{2,2}^d - x_{12,}^d \end{matrix}\right)$$
of ${\tt T}$, and we introduce variables $(\begin{smallmatrix} a_1 \\ a_2\end{smallmatrix})  = {\tt T}_1 \bm \alpha_1  + {\tt T}_2 \bm \alpha_2  \in \Bbb{R}^2$ with $(\begin{smallmatrix} \bm \alpha_1 \\ \bm \alpha_2\end{smallmatrix}) \in \Bbb{R}^{2 + (k-2)}.$ The integral over $\bm \alpha_2$ is $O(1)$, the integral over $a_2 =  q_{\bm \alpha}(\bm{x}_2) - q_{\bm \alpha}(\bm{x}_1) \ll M^{d-k}$ is $O(M^{d-k})$, and by Cramer's rule the integral over $a_1$ is $\ll |\det{\tt T}_1 |/ |x_{2,1}^d - x_{1,1}^d| \ll |\det {\tt T}_1|/ M^{d-1}X$. Thus we have in total  
$\mathcal{I}^{(2)}(\bm{x}, \bm{y}, M) \ll M^{1-k} /X$. 

Suppose first that $\rho_2 = \sigma_2 = 0$. Then we determine the first column trivally in $O(M^3X)$ ways which determines $\rho_1$ and $\sigma_2$ so that every other column can be chosen in $O(M^2)$ ways. This gives a  contribution
$$ \frac{M^{1-k}}{X} M^3 X M^{2(k-1)} =  M^{k+2}.$$

After possible relabelling, we assume from now on without loss of generality $\rho_1\rho_2 \not = 0$. 
 Choose the first two columns in $M^7X$ ways. 
 This determines $\rho_1, \rho_2, \sigma_1, \sigma_2$ since the upper $2$-by-$2$ block  is invertible. 
For every other column we apply Lemma \ref{lemma-a}a with $k=3$ to the first equation in \eqref{rhosigma} to determine $x_{1,i}, x_{2,i}, y_{1,i}$, and then $y_{2,i}$ is determined from the second equation in \eqref{rhosigma}. This gives a contribution
$$ \frac{M^{1-k}}{X} M^7X M^{(1+\varepsilon)(k-2)} =  M^{6+\varepsilon}$$
and completes the proof.

\begin{lemma}\label{c3rank3}
	For $d\geq 2, k\geq \max(d, 4)$ we have 
	$$ \mathcal C_{3}^{(2)}(M; k,2, d)\ll M^{\varepsilon}(M^{2k-1}+M^{2d+2})$$
	for any $\varepsilon >0$.
\end{lemma}
\textbf{Proof.} We may assume without loss of generality that 
\begin{equation}\label{1}
y_{1,i}^d = \rho_1 x_{1,i}^d + \rho_2 x_{2,i}^d + \rho_3 y_{2,i}^d
\end{equation}
for some $\rho_i\in \mathbb{R}$.
As a change of variable matrix we can take the submatrix with entries in the first three columns and $1, 2,4$-th rows of ${\tt T}$, i.e.\ 
\begin{align}\label{T1rank3}
{\tt T}_1=   \left(\begin{matrix} x_{1,1}^d & \ldots & x_{1,3}^d\\ x_{2,1}^d - x_{1,1}^d & \ldots & x_{2,3}^d - x_{1,3}^d\\  y_{2,1}^d - y_{1,1}^d & \ldots & y_{2,3}^d - y_{1,3}^d \end{matrix}\right)
\end{align}
whose determinant is of size $D_3$ by \eqref{XYDef}. 
We first consider the case when $D_2=0$. If $D_2=0$, then we have
\begin{align}\label{D2=0}
	x_{2,k}^d-x_{1,k}^d=\rho(y_{2,k}^d-y_{1,k}^d), \text{ for $k=1, 2,3$ for some $\rho\not=0$.}
\end{align}  We can choose the first colum in $O(M^3X)$ ways and they determine $\rho$. For all other columns, we see that \eqref{D2=0} is a diagonal form in $4$ variables and Lemma \ref{lemma-a} gives $O(M^{2+\varepsilon})$ choices uniformly for all $\rho$. We have that $\mathcal I^{(2)}(\bm x,\bm y, M)\ll M^{1-k}/X$ by using the matrix $\mathtt{T}_1$ in Lemma \ref{c2rank2}. Thus we have a total contribution 
$$\frac{M^{1-k}}{X}M^3X M^{2(k-1)+\varepsilon}=M^{k+2+\varepsilon}\ll M^{2k-1+\varepsilon},$$
using $k\geq 4$.  

Now we consider the case when $D_2\not=0$. 
Since $\rank {\tt T}=3$, we have $D_3>0$. As before we see (using Cramer's rule) that  
$$\mathcal{I}^{(2)}(\bm{x}, \bm{y}, M) \ll M^{2d-2k}/D_2.$$

Let us first deal with the case $Y = 0$.  The first column has $M^2X$ choices since $y_{2,1}=y_{1,1}$. In the second column we pick $x_{1,2}, x_{2,2}, y_{1,2}$ randomly. The second condition in \eqref{oa2}  now reads
$$|(x_{2,1}^d - x_{1,1}^d)(y_{2,2}^d - y_{1,2}^d)| \asymp D_2.$$
Thus $0 \not= |y_{2,2} - y_{1,2}| \ll D_2/XM^{2d-2}$. Treating all other columns trivially by $O(M^4)$, we obtain a contribution
\begin{equation}\label{trivial}
\frac{M^{2d-2k}}{D_2}   \cdot M^2 X  \cdot M^3 \frac{D_2}{XM^{2d-2}}  \cdot M^{4(k-2)} = M^{2k-1}.
\end{equation}
(This argument could be improved in the present case, but we will reuse it in the next lemma.) 

From now on assume $Y > 0$. We first consider the case when $\rho_2 = \rho_3 = 0$ (or similarly $\rho_1 = \rho_3 = 0$). Choose the first column in $O(M^2XY)$ ways using \eqref{XYDef}. This determines $\rho_1\not= 0$ and so that we have trivially $O(M^3)$ choices for all other columns.  Alternatively, we can forget about $\rho_1$ and  use the second condition in \eqref{oa2} in the form
\begin{equation}\label{deteq}
|(x_{2,1}^d - x_{1,1}^d)(y_{2,i}^d -  y_{1,i}^d) - (y_{2,1}^d - y_{1,1}^d)(x_{2,i}^d - x_{1,i}^d) | \ll D_2.
\end{equation}
Using Lemma \ref{lemma-b}b together with the trivial bound, the number of choices for the $i$-th column can be bounded by 
\begin{align}
  M^{\varepsilon}\min \Big(\frac{D_2M^{4-d}}{M^{d-1} \sqrt{XY}} + \frac{D_2^{1/2} M^{3-d/2}}{M^{(d-1)/2} Y^{1/2}} + M^2 , M^3\Big).
\end{align}
This gives a contribution towards $\mathcal C_3^{(2)}(M)$ when $k\geq 3$
\begin{equation}\label{small}
\begin{split}
&\frac{M^{2d - 2k+\varepsilon}}{D_2} M^2XY  \Big( \min\Big(\frac{D_2M^{4-d}}{M^{d-1} \sqrt{XY}} + \frac{D_2^{1/2} M^{3-d/2}}{M^{(d-1)/2} Y^{1/2}}, M^3\Big) + M^2 \Big)^{k-1}\\
& \ll \frac{M^{2d-2k+\varepsilon}}{D_2}M^2 XY \Big( \frac{D_2M^{4-d}}{M^{d-1}\sqrt{XY}}M^{3(k-2)}+\frac{D _2 M^{6-d}}{M^{d-1}Y}M^{3(k-3)}+M^{2(k-1)}\Big)\\
& \ll M^{\varepsilon}\Big(M^{k+2}+\frac{M^{2d+2}}{D_2}\Big).
\end{split}
\end{equation}

From now on we can assume that \eqref{1} is at least a ternary form (i.e.\ at least two $\rho_i$'s are non-zero).  We fix the first column randomly in $O(M^2XY)$ ways. For the next two columns we apply Lemma \ref{lemma-b}b to the determinant equation \eqref{deteq}.  Then $\rho_1, \rho_2, \rho_3$ are determined since the relevant $3 \times 3$ block  ${\tt T}_1$ of $\tt T$ in \eqref{T1rank3} is invertible. For the rest of the columns we can apply Lemma \ref{lemma-a}a using \eqref{1} with $k=3,4$ depending on whether $\rho_1\rho_2\rho_3=0$ or not. This gives a total contribution
\begin{equation}\label{big} 
\begin{split}
& M^\varepsilon \frac{M^{2d-2k}}{D_2} M^2XY \Big( \frac{D_2M^{4-d}}{M^{d-1} \sqrt{XY}} + \frac{D_2^{1/2} M^{3-d/2}}{M^{(d-1)/2} Y^{1/2}} + M^2 \Big)^2 M^{2(k-3)}\\
&  \ll M^{\varepsilon}\Big(M^{6 - 2d} D_2 + XM^3 + \frac{M^{2d}}{D_2} XY\Big) \ll  M^{2d + 2+\varepsilon}.
\end{split}
\end{equation}
Combining \eqref{small} and \eqref{big} together with \eqref{trivial} completes the proof. 

 \subsection{Full rank}
 Now we consider the case $\rank{\tt T}=4$ and so we can assume $D_2, D_3, D_4>0$ in \eqref{XYDef} from assumptions in \eqref{oa2}.
 
 \begin{lemma}\label{c4rank4}
	For $d \geq 2, k \geq \max(4, d)$ we have 
	$$ \mathcal C_{4,-}^{(2)}(M, D_0)\ll M^{\varepsilon}\Big(M^{2k} \frac{D_0^{1/7}}{M^{4d/7}} + M^{2k-1} + M^{(7k+3)/4} + M^{2d+2}\Big)$$
	for any $\varepsilon >0$.
\end{lemma}

\textbf{Proof.} In the case $\rank({\tt T}) = 4$, the usual change of variables in \eqref{quantity2} is 
$$(a_1, a_2, a_3, a_4)^{\top} = {\tt T}\bm \alpha = (q_{\bm \alpha}(\bm{x}_1), q_{\bm \alpha}(\bm {x}_2)-q_{\bm \alpha}(\bm{x}_1), q_{\bm \alpha}(\bm{y}_1), q_{\bm \alpha}(\bm{y}_2) - q_{\bm \alpha}(\bm{y}_1))^{\top}$$
and $\bm a  = {\tt T}_1 \bm \alpha_1  + {\tt T}_2 \bm \alpha_2$ with the notation as in \eqref{Txy} and the subsequent line. For the integral over $a_1, a_3$ we use that the support $F$ implies that ${\tt T}_1^{-1}(a_1, *, a_3, *)^{\top}$ lies in a box of length $O(1)$ (depending on $a_2, a_4$), and by a (generalized version of) Cramer's rule, this volume is at most $\ll |\det {\tt T_1}|/D_2$.  
The integral over $\bm \alpha_2$ is $O(1)$ and the integral over $a_2, a_4$ is $O(M^{2d-2k})$. so that in total   $$ \mathcal{I}^{(2)}(\bm{x}, \bm{y}, M) \ll \frac{M^{2d}}{D_2} M^{-2k}.$$

The case $Y=0$ yields verbatim as in \eqref{trivial} a contribution $O(M^{2k-1})$, and 
from now on we   assume $Y > 0$. We have various bounds depending on the size of $D_2$.

For the $i$-th column we  consider \eqref{deteq} as an inequality in $x_{1, i}, x_{2, i}, y_{1, i}, y_{2, i}$. 
We can apply Lemma \ref{lemma-b}b for the $i$-th column or alternatively the trivial bound $O(M^4)$ getting
\begin{equation}\label{2}
\begin{split}
&M^\varepsilon \frac{M^{2d - 2k}}{D_2} M^{2}XY \Big(\min\Big( \frac{D_2M^{4-d}}{M^{d-1} \sqrt{XY}} + \frac{D_2^{1/2} M^{3-d/2}}{M^{(d-1)/2} Y^{1/2}}, M^4\Big) + M^2 \Big)^{k-1}\\
& \ll M^\varepsilon \frac{M^{2d-2k}}{D_2}M^2 XY \Big(\Big( \frac{D_2^2M^{2(4-d)}}{M^{2(d-1)}XY}+\frac{D_2M^{6-d}}{M^{d-1}Y}\Big)M^{4(k-3)}+M^{2(k-1)}\Big)\\
& \ll M^\varepsilon \Big(\frac{M^{2k} D_2}{M^{2d}} + M^{2k-2} + \frac{M^{2d+2}}{D_2}\Big).
\end{split}
\end{equation}
This is useful when $D_2\ll M^{2d-\delta}$ for some $\delta > 0$. 

When $D_2$ is large, we use information from $3 \times 3$ determinants as described in the third condition in \eqref{oa2}.  Here we  choose the first two columns randomly in $O(M^8)$ ways. Then for any other column we apply Lemma \ref{lemma-b}a with $k=4$ getting $O(M^{4-d} D_3/D_2 + M^3)$ choices using the third condition in \eqref{oa2}. Indeed, the coefficient of $x_{1, i}^d$ is $\check{\mathcal{T}}^{\{1, 2\}}_{\{2, 4\}} + \check{\mathcal{T}}^{\{1, 2\}}_{\{1, 4\}}$ and the coefficient of $x_{2, i}^d$ is $-\check{\mathcal{T}}^{\{1, 2\}}_{\{1, 4\}}$, so the maximum of the absolute values of the coefficients is $\gg D_2$. Alternatively we can chose the $i$-th column trivially in $M^4$ ways. When $k \geq 3$ this gives a contribution  
\begin{displaymath}
\begin{split}
& \frac{M^{2d-2k}}{D_2} M^8 \Big(\min\Big( \frac{M^{4-d} D_3}{D_2}, M^4\Big) + M^3\Big)^{k-2}\\
& \ll  \frac{M^{2d-2k}}{D_2} M^8 \Big(  \frac{M^{4-d} D_3}{D_2} M^{4(k-3)}  + M^{3(k-2)}\Big) \ll \frac{M^{2k+d} D_3}{D_2^2} +  \frac{M^{k + 2d + 2}}{D_2}.
\end{split}
\end{displaymath}
Combining this with \eqref{2}, we obtain for $k\geq 3$ that 
\begin{equation}\label{3}
 M^{\varepsilon}\big(M^{2k-d} D_3^{1/3} +  M^{1 + 3k/2} + M^{2k-2}+M^{2d+2}\big).
 \end{equation}
This is useful when $D_3\ll M^{3d-\delta}$. 
When $D_3$ is large, we use instead the fourth condition in \eqref{oa2}, in the exact same way we also have the bound (observing that $D_3 \ll D_2 M^d$ by \eqref{oa2})
\begin{displaymath}
\begin{split}
& \frac{M^{2d-2k + d}}{D_3} M^{12} \Big(\min\Big( \frac{M^{4-d} D_4}{D_3}, M^4\Big) + M^3\Big)^{k-3}\\
& \ll  \frac{M^{2d-2k+d}}{D_3} M^{12} \Big(  \frac{M^{4-d} D_4}{D_3} M^{4(k-4)}  + M^{3(k-3)}\Big) \ll \frac{M^{2k+2d} D_4}{D_3^2} +  \frac{M^{k + 3d + 3}}{D_3}.
\end{split}
\end{displaymath}

Combining this with \eqref{3} gives the final bound
\begin{equation*} 
M^{\varepsilon}\Big(M^{2k} \frac{D_4^{1/7}}{M^{4d/7}}  + M^{(7k+3)/4} +M^{2k-2}+ M^{2d+2}\Big). 
 \end{equation*}
 
Recalling \eqref{trivial} for the case $Y=0$, we complete the proof of the lemma.  
  
  \subsection{The linear case}
  In order to complete the proof of Theorem \ref{thm2} based on Proposition \ref{prop3}, we also need to verify \eqref{prop2-bound} for $\ell=2, k\geq d$ and some $0<\eta, \delta'<1$. This is much simpler than the analysis in Lemmas \ref{c2rank2} -- \ref{c4rank4} and can be dealt with quickly. 
  \begin{lemma}
  	For $k \geq d\geq 2 = \ell$, we have
  	\begin{align}
  	\mathcal C_{2,-}(M, D_0)\ll M^{k+\varepsilon} \frac{D_0}{M^{2d}}
  	\end{align}
  	for any $\varepsilon>0$.
  \end{lemma}
  \textbf{Proof.}
  Here we define $T$ as in \eqref{TT'} with $\ell=2$. Since we must have $\rank T=2$, we can rearrange columns of $T$ so that
  $X\asymp |x_{2,1}-x_{1,1}|\geq |x_{2,i}-x_{1,i}|$ for all $i\geq 2$. Using the change of variable matrix $T_1$ defined by
  $$T_1=\begin{pmatrix}
  x_{1,1}^d & x_{1,2}^d\\
  x_{2,1}^d-x_{1,1}^d & x_{2,2}^d-x_{1,2}^d
  \end{pmatrix}$$
  we see that 
  $$\mathcal{I}(\bm{x}, M) \ll  \frac{1}{|\det T_1|}\frac{M^{d-k}|\det T_1|}{M^{d-1} X}=\frac{1}{M^{k-1} X}.$$
  Note that $|\det T_1 |\ll D_0$ implies $$ \frac{D_0}{M^{2d-2}} \gg |x_{1,1} x_{2,2} - x_{1,2} x_{2,1}| =| x_{1,1}( x_{2,2} - x_{1,2}) - x_{1,2} (x_{2,1} - x_{1,1})|. $$
  We fix $|x_{2,2} - x_{1,2} |\ll X$ and $|x_{1,1}| \ll M$ and any integer $\ll  D_0/M^{2d-2}$, then by a divisor argument $x_{1,2}\not= 0$ and $x_{2,1} - x_{1,1} \not= 0$ are fixed up to a factor $M^{\varepsilon}$. Hence the first two columns can be chosen in 
     $\ll M^\varepsilon MX  D_0/M^{2d-2}$ ways.  All other columns can be bounded trivially in $O(M^2)$ ways, so that   \begin{align}
  \mathcal C_{2,-}(M, D_0)&\ll \frac{M^\varepsilon}{M^{k-1}X}MX \frac{D_0}{M^{2d-2}}M^{2(k-2)}\ll M^{k+\varepsilon} \frac{D_0}{M^{2d}}\   \end{align}
  as desired. 
  
 \subsection{The case $d=2$} Finally we prove the second half of Theorem \ref{thm5}, so we assume $d = \ell = 2$. The results in the previous sections, in particular Lemmas  \ref{c2rank2} -- \ref{c4rank4} suffice for the case $k\geq 4$, so from now on we assume $k=3$.  This case is notationally different from the general set-up in Section \ref{sec3} in that the definition of $\tt T$ in \eqref{Txy} requires $k \geq 2\ell = 4$ for the change of variable matrix to exists. Here we consider instead the matrix
 $$\tt T = \left(\begin{matrix}  x_{2,1}^2 - x_{1,1}^2 & x_{2,2}^2-x_{1,2}^2 & x_{2,3}^2 - x_{1,3}^2
 \\ y_{2,1}^2 - y_{1,1}^2 & y_{2,2}^2-y_{1,2}^2 & y_{2,3}^2 - y_{1,3}^2 \end{matrix}\right).$$
 We can arrange the columns so that 
 \begin{equation}\label{d2Toa}
 \begin{split}
& |x_{2,1}-x_{1,1}|\geq |x_{2,i}-x_{1,i}|,\\
 &X \asymp  |x_{2,1}-x_{1,1}|\geq |y_{2,1}-y_{1,1}|  \asymp Y,\\
  &D_2 \asymp |\dot{\mathcal{T}}^{\{1, 2\}}_{\{1, 2\}} | \geq |\dot{\mathcal{T}}^{\{1, i\}}_{\{1, 2\}} |, \ i\geq 2.
 \end{split}
 \end{equation}
 where $\dot {\mathcal T}^{\mathfrak c}_{\mathfrak r}$ is the determinant of the submatrix $\tt T$ with entries in columns in $\mathfrak c$ and rows in $\mathfrak r$. 
 
Let us first assume that $\text{rank}({\tt T}) = 2$ so that $D_2>0$. After the change of variables $$\textbf{a} = (a_1,a_2)^{\top} =(q_{\bm \alpha}(x_2)-q_{\bm \alpha}(\bm x_1), q_{\bm \alpha}(y_2)-q_{\bm \alpha}(\bm y_1))^{\top} ,$$ the integration over $a_1, a_2$  in  $\mathcal{I}^{(2)}(\bm{x}, \bm{y}, M)$ defined in   \eqref{quantity2}
is $O(M^{2(d-k)})$, so that the entire integral is bounded by 
 $$\mathcal{I}^{(2)}(\bm{x}, \bm{y}, M) \ll M^{2(d-k)}D_2^{-1}=(M^2 D_2)^{-1}.$$ 
 
 Let us quickly deal with the case $Y=0$ as in \eqref{trivial}. Here we choose the first column in $O(M^2X)$ ways and we choose the remaining four  $x$-variables in $O(M^2X^2)$ ways by the first condition in \eqref{d2Toa}. The third condition in \eqref{d2Toa} now implies that $|y_{2, i} - y_{1, i}|\ll D_2/M^2 X$ for $i = 2, 3$. Note that $y_{2, 2} - y_{1, 2} \not = 0$ since otherwise the conditions $Y=0$ and $q_{\bm \alpha}(\bm y_2)-q_{\bm \alpha}(\bm y_1)\ll M^{d-k}$ imply  $\bm y_1=\bm y_2$. Thus we get a total contribution of
 \begin{equation}\label{y0}
 \frac{1}{M^2 D_2} M^2X \cdot M^2X^2  \cdot M^2 \frac{D_2}{M^2X} \Big(1+\frac{D_2}{M^2X}\Big) = M^2X^2 + D_2X \ll M^4 + D_2M.
 \end{equation}
    
From now on we assume $Y\not= 0$. The first column can be chosen in $O(M^2XY)$ ways.  Applying Lemma \ref{lemma-b}b with $d=2$ to the third condition in \eqref{d2Toa}, we can bound the number of choices for the second and third column by 
 \begin{align}\label{d2XY}
 M^\varepsilon\Big(\frac{M^2 D_2}{MX}+\frac{MD_2}{MY}+M^2\Big).
 \end{align}
 Therefore   a bound towards $\mathcal C_2^{(2)}(M)$ satisfying  \eqref{d2Toa} is given by 
 \begin{align}
 \frac{M^{\varepsilon}}{M^2 D_2}M^2XY\Big( \frac{M^2 D_2}{M X}+ \frac{MD_2}{MY}+M^2\Big)^2
 &\ll M^{\varepsilon}\Big( M^2 D_2+ D_2+\frac{M^6 }{D_2}\Big).
 \end{align}
 which together with \eqref{y0} gives (using $D_2\ll M^4$)
 \begin{align}\label{middleT2}
 M^\varepsilon (M^2 D_2+ M^4+ \frac{M^6}{D_2}).
 \end{align}

 We complement this with the following alternative bound under the assumption $D_2 \ll M^2$ with a sufficiently small constant. This plays the same role as the condition \eqref{D2bound}, and as in that discussion we see that in this case none of the entries in the first two columns of ${\tt T}$ can vanish as $D_2\not= 0$.  Fix $g = (x_{2,1}^2 - x_{1,1}^2, y_{2,1}^2 - y_{1,1}^2) \ll MY$. Then the first column of ${\tt T}$ can be chosen in $O(M^2XY/g^2)$ ways. For the second column we consider the determinant equation
 $$|a (y_{2,2}^2 - y_{1,2}^2) - b (x_{2,2}^2 - x_{1,2}^2) | \asymp D_2/g, \quad a = (x_{2,1}^2 - x_{1,1}^2)/g, \, b = (y_{2,1}^2 - y_{1,1}^2)/g.$$
 Since $(a, b) = 1$, by a divisor argument we find at most $D_2M^{2+\varepsilon}/g|a|$ solutions for the second column. Using \eqref{d2XY} for the third column  we obtain a contribution in total
 $$\frac{M^{\varepsilon}}{M^2 D_2}\sum_{g \ll MY} \frac{M^2XY}{g^2} \Big( \frac{D_2 M^2}{MX} \Big)\Big(\frac{M^2 D_2}{MX}+\frac{MD_2}{MY}+M^2\Big) \ll  M^{\varepsilon}\Big(M^{2} D_2+M^{4}\Big).$$
 We choose this bound if $D_2 \ll M$ (in which case it is applicable) and \eqref{middleT2} otherwise, obtaining in total 
 $$\mathcal C_{2,-}^{(2)}(M, D_0) \ll M^{6+\varepsilon}\Big( \frac{1}{M} + \frac{D_0}{M^4}\Big)$$
 which is admissible.
 
 It remains to deal with the case $\text{rank}({\tt T}) = 1$.  The integral in \eqref{quantity2} can be bounded by the same argument by $\ll 1/M^2X$.  Now we simply choose the first column in $O(M^3X)$ ways. Since the two rows of ${\tt T}$ are linearly dependent, this determines the factor of proportionality. Note that if $y_{2,1}=y_{1,1}$, then we will have $\bm y_{2}=\bm y_1$ which is forbidden.  Thus we can apply Lemma \ref{lemma-a}a to see every other column can then be chosen in $O(M^{2+\varepsilon})$ ways, and we obtain a total contribution  
 $$\mathcal C_1^{(2)}(M) \ll  \frac{  M^3X \cdot M^{4+\varepsilon}}{M^2X} =M^{5+\varepsilon}$$
 which is admissible.

\end{document}